\documentclass[VANCOUVER,STIX1COL]{WileyNJD-v2}

\usepackage{pgfplots}
\usepgfplotslibrary{fillbetween}
\pgfplotsset{compat=newest,
every axis/.append style={axis x line=bottom,
                          axis y line=left,
                          scale only axis,
                          y label style={at={(0.0,1.0)},anchor=south west,rotate=-90}
                          },
}
\pgfrealjobname{CVaR_ML_Rev4}

\newcommand{\be}{\begin{eqnarray}}
\newcommand{\e}{\end{eqnarray}}
\newcommand{\grad}[2]{\nabla_{#2}#1}

\articletype{Article Type}%

\date{1 Dec 2022}

\begin{document}
\title{TTRISK: Tensor Train Decomposition Algorithm for Risk Averse Optimization}
\author[1]{Harbir Antil}
\address[1]{\orgdiv{The Center for Mathematics and Artificial Intelligence
(CMAI) and Department of Mathematical Sciences}, \orgname{George Mason University},
\orgaddress{Fairfax, VA 22030, \country{USA}}. \email{hantil@gmu.edu}}

\author[2]{Sergey Dolgov}
\address[2]{\orgdiv{Department of Mathematical Sciences}, \orgname{University of Bath}, \orgaddress{Bath, BA2 7AY, \country{UK}}. \email{s.dolgov@bath.ac.uk}}

\author[3]{Akwum Onwunta}
\address[3]{\orgdiv{Department of Industrial and Systems Engineering}, \orgname{Lehigh University}, \orgaddress{Bethlehem, PA 18015, \country{USA}}. \email{ako221@lehigh.edu}}

\authormark{Harbir Antil \textsc{et al}}

\corres{*Sergey Dolgov, University of Bath, Bath, BA2 7AY, UK. \email{s.dolgov@bath.ac.uk}}


\abstract[Summary]{This article develops a new algorithm named TTRISK to solve high-dimensional
risk-averse optimization problems governed by differential equations (ODEs and/or
PDEs) under uncertainty. As an example, we focus on the so-called Conditional 
Value at Risk (CVaR), but the approach is equally applicable to other coherent 
risk measures. Both the full and reduced space formulations are considered.
The algorithm is based on low rank tensor approximations of random
fields discretized using stochastic collocation.
To avoid non-smoothness of the objective function underpinning the CVaR,
we propose an adaptive strategy to select the width parameter of the
smoothed CVaR to balance the smoothing and tensor approximation errors.
Moreover, unbiased Monte Carlo CVaR estimate can be computed by using
the smoothed CVaR as a control variate. To accelerate the computations, 
we introduce an efficient preconditioner for the Karush-Kuhn-Tucker (KKT) system in the
full space formulation.The numerical experiments demonstrate that
the proposed method enables accurate CVaR optimization constrained 
by large-scale discretized systems. In particular, the first example consists 
of an elliptic PDE with random coefficients as constraints. The second example 
is motivated by a realistic application to devise a lockdown plan for United 
Kingdom under COVID-19. The results indicate that the risk-averse framework is feasible with the tensor approximations under tens of random variables.}


\jnlcitation{\cname{%
\author{Harbir Antil},
\author{Sergey Dolgov}, and
\author{Akwum Onwunta}} (\cyear{2022}),
\ctitle{TTRISK: Tensor Train Decomposition Algorithm for Risk Averse Optimization}, \cjournal{NLAA}, \cvol{2022}.}

\maketitle

\section{Introduction}

Uncertainty is ubiquitous in science and engineering applications. It may arise due to noisy 
measurements, unknown parameters, or unverifiable modeling assumptions. Examples 
include, infectious disease models for COVID-19  \cite{DGKP-SEIR-2021}, partial differential 
equations (PDEs) with random coefficients, boundary conditions or right-hand-sides 
\cite{durlofsky2012uncertainty,petrat2012inexact,DMTartakovsky_AGuadagnini_MRiva_2003a,DPKouri_TMSurowiec_2020a}.
The control or design optimization problems constrained by such systems must produce controls 
or optimal designs which are resilient to this uncertainty.  To tackle this, recently in \cite{KS16,DPK2018,DPKouri_TMSurowiec_2020a,GSU21}, the authors have created  
risk-averse optimization frameworks targeting engineering applications. 

The goal of this paper is to introduce a new algorithm TTRISK which uses a Tensor Train (TT) 
decomposition to solve risk averse optimization problems constrained by differential equations 
(ODEs and/or PDEs). Let $(\Omega, \mathcal{A},\mathbb{P})$ be a complete probability space.
Let $U, Y$ be real reflexive Banach spaces, and let $Z$ be a real Banach
space. Here $Y$ denotes the deterministic state space, $U$ is the space of optimization variables
(control or designs etc.) and $Z$ is the differential equation residual space. Let $U_{ad} \subseteq
U$ be a closed convex subset and let $c : Y \times U_{ad} \times \Omega \rightarrow Z$ denote, e.g.,
a partial differential operator, then consider the equality constraint
	\[
		c(y,u;\omega)=0, \quad \mbox{in } Z, \quad \mbox{a.a. }  \omega \in \Omega 	
	\]
where a.a. indicates ``almost all" with respect to a probability measure $\mathbb{P}$.
The goal of this article is to consider optimization problems of the form
\begin{equation}
\label{modprob_intro}
	\min_{u\in U_{ad}}\mathcal{R}[\mathcal{J}(y,u;\omega)]
		+ \alpha \mathcal{P}(u) \;\; \mbox{subject\; to}\;\;  c(y,u;\omega)=0,
		\quad \mbox{in } Z, \quad \mbox{a.a. }  \omega \in \Omega,
\end{equation}
where $u \in U_{ad}$ is the deterministic control and $y \in Y$ is the state, $\mathcal{P}$ 
is the cost of the control, $\alpha \ge 0$ is the regularization parameter, $\mathcal{J}$ is the 
uncertain variable objective function and $\mathcal{R}$ is 
the risk-measure functional which maps random variables to extended real numbers. 

We assume that $\mathcal{R}$ is based on expectation, i.e., 
	\begin{equation}
	\label{eq:risk}
		\mathcal{R}[X]=\inf_{t\in \mathcal{T}}\mathcal{R}_t[X], \quad \mbox{where} \quad \mathcal{R}_t[X]:=\mathbb{E}[f(X,t)],
	\end{equation}
$f:\mathbb{R}\times \mathbb{R}^N\rightarrow \mathbb{R}$  and
$\mathcal{T} \subseteq\mathbb{R}^N$, with $N\in\mathbb{N}$, is a closed convex set.
The problem class \eqref{eq:risk} consists of a large number of risk-measures
that are of practical interest.
In particular, it includes the \emph{coherent risk measures}, 
which are sub-additive, monotonic, translation equivariant and positive 
homogeneous \cite{RU2000,AShapiro_DDentcheva_ARuszczynski_2014a}. 
Notice that sub-additivity and positive homogeneity implies convexity. 
These risk measures have several advantages, for instance, they preserve desirable 
properties of the original objective function such as convexity. In addition, in engineering 
or in finance applications, tail-probability events may be rare but critical if they lead to failure of a system. It is therefore essential to  minimize the risk of
failure, i.e. $\mathcal{R}$, and obtain controls $u$ which are resilient to uncertainty
in the system.

A typical example of coherent risk measure $\mathcal{R}$ is the conditional value-at-risk 
($\mbox{CVaR}_{\beta}$), where $f$ in \eqref{eq:risk} is given by
\begin{equation}\label{eq:fcvar}
	f(X,t)= t+(1-\beta)^{-1}(X-t)_+ ,
\end{equation}
with $\mathcal{T}=\mathbb{R}$, $\beta \in (0,1)$ is the confidence level and $(x)_+=\max\{x,0\}$.  
CVaR$_{\beta}$ is also known as expected shortfall. It's origin lies in financial 
mathematics \cite{RU2000, RU2002}, but owing to Kouri \cite{kouri2012approach} and Kouri
and Surowiec \cite{KS16}, it is now being widely used in engineering applications. 
Our work in particular, focuses on minimization problems \eqref{modprob_intro} with $\mathcal{R}$
given by CVaR$_{\beta}$ but it can be extended to other coherent risk measures, such as 
buffered probability of exceedence (BPOE), of type \eqref{eq:risk}.

Notice that, since risk measures, such as CVaR$_{\beta}$ focus on the upper tail events, the 
traditional sampling techniques to solve these stochastic PDE-constrained optimization problems are 
often computationally expensive. More precisely, CVaR$_{\beta}$ captures the cost associated with 
rare events, but it requires more samples in order to be accurately approximated, which leads to
many differential equation solves \cite{KS16}. 
Moreover, the presence of the non-smooth function $(\cdot)_+$ in CVaR$_{\beta}$ poses several 
challenges, including, nondifferentiable cost functional, wasted Monte Carlo samples outside of the 
support of $(\cdot)_+ = \max\{\cdot, 0 \}$, or slowly converging polynomial and other function 
approximation methods.

To tackle some of these challenges, \cite{KS16} has proposed a smoothing of 
$(\cdot)_+$ which requires solving a sequence of smoothed optimization problems 
using Newton-based methods. Another solution strategy is to reformulate the problem 
and use interior-point methods \cite{GSU21}. A duality based approach has also been 
recently proposed in \cite{kouri2021primal}.

In this paper we develop an efficient method to tackle the above challenges associated 
with minimization of CVaR$_{\beta}$ subject to constraints given by differential equations 
with random inputs.
We consider two formulations
of \eqref{modprob_intro}. The first one is the implicit approach where we remove the
equality constraint $c(y,u;\omega) = 0$ via a control to solution map $u \mapsto y$ 
\cite{KS16,GSU21,kouri2021primal}. The second case is the full space approach, where
we directly tackle the full problem \eqref{modprob_intro} using the Lagrangian
formulation. The latter
formulation appears to be new in the context of risk-averse optimization.
Numerical experiments demonstrate that the full formulation converges more reliably for extreme parameters, e.g. large $\beta$ and small $\alpha$.

Our framework builds rigorously on tensor decomposition methods, which emerged in the past two decades \cite{hackbusch-2012,khor-book-2018} as an efficient approximation of multi-index arrays, in particular when those contain expansion coefficients of high-dimensional functions \cite{Marzouk-stt-2016,Gorodetsky-ctt-2019}.
The idea starts from the classical separation of variables.
Functions of certain structure \cite{khor-low-rank-kron-P1-2006} or regularity \cite{uschmajew-approx-rate-2013} have been shown to admit rapidly (often exponentially) convergent series, where each term is a product of univariate functions.
Later, instead of a simple sum of products,
it was found more practical to consider hierarchical separation of variables \cite{hackbusch-2012}.
A particularly simple instance of such is the Tensor Train (TT) decomposition \cite{osel-tt-2011} that admits efficient numerical computations.
One of the most powerful algorithms of this kind is the cross approximation, as well as its variants \cite{ot-ttcross-2010,so-dmrgi-2011proc,ds-parcross-2020}.
Those allow one to compute a TT approximation to potentially any function, using a number of samples from the sought function that is a small multiple of the number of degrees of freedom in the tensor decomposition.
Once TT decompositions are computed, integration, differentiation and linear algebra of the original functions can be implemented using their TT formats instead with a linear cost in the dimension.

However, irregular functions, such as the $(\cdot)_+$ function in \eqref{eq:fcvar}, may lack an efficient TT decomposition.
The main novelty of the paper is an algorithm that is adaptive in both the TT complexity and width parameter in the smoothed CVaR function,
which allows one to actually alleviate the curse of dimensionality,
since smooth functions do admit convergent TT approximation.
If the bias from the smoothing is still too large for a feasible TT decomposition,
we can obtain an unbiased, asymptotically exact solution with a version of
Multilevel Monte Carlo methods \cite{giles-mlmc-2015}, namely, we use a smoothed solution as a control variate \cite{robert-MC-book-2004}.

The numerical experiments demonstrate that the stochastic risk-averse control problem can be solved 
with a cost that depends at most polynomially on the dimension. This allows us to solve a realistic 
risk-averse ODE control problem with $20$ random variables.

\medskip
\noindent 
{\bf Outline:} Section~\ref{s:back} and Appendix set up the relevant notation and provide the necessary
background on risk-averse-optimization and tensor-train decomposition.
In Section~\ref{sec:reduced},
we introduce the
control-to-state map, i.e., $u \mapsto y$, and eliminate the equality constraints $c(\cdot) = 0$.
The resulting optimization problem \eqref{modprob_intro} is only a function of
the control variable $u$ and is known as reduced problem.
A control variate correction of the problem is considered in Section~\ref{s:control_variate}.
When the constraint $c(\cdot)$ is handled directly, the resulting formulation is called full-space problem
and is discussed in Section~\ref{sec:lagrangian}. This is followed by Section~\ref{s:GN} where
the Gauss-Newton system for the full-space formulation is considered. Next, in 
Section~\ref{s:precond} we consider preconditioning strategies for this formulation.
Section~\ref{s:numerics} 
focuses on our numerical examples, where we first consider an optimal control problem 
constrained by an elliptic PDE with random coefficients. This is followed by the risk-averse 
optimal control of an infectious disease model which has been recently developed to propose lockdown strategies in the United Kingdom due to COVID-19.

\section{Background}
\label{s:back}

In this section, we first provide background on CVaR$_\beta$ and introduce a regularized problem 
with CVaR$_\beta$ replaced by CVaR$_\beta^\varepsilon$ with $\varepsilon > 0$. This is followed by a 
discussion on tensor-train decomposition. The aim of this section is to set the stage for TTRISK.

\subsection{Risk-averse optimization}

Let $\mathcal{L}:=(\Omega,\mathcal{A},{\mathbb{P}})$ be  a complete probability space.
Here, $\Omega$ is the set of outcomes, $\mathcal{A}\subset 2^{\Omega}$ is the $\sigma$-algebra of events, and 
$\mathbb{P}:\mathcal{A}\rightarrow [0,1]$ is an appropriate probability measure. Let $X$ be a (scalar) random variable
defined on $\mathcal{L}.$
For example, we will consider the objective function $\mathcal{J}$ in what follows.
Then, the expectation of $X$ denoted by $\mathbb{E}[X]$ is given by
\[
\mathbb{E}[X]=\int_{\Omega}\; X(\omega) d\mathbb{P}(\omega). 
\]
As stated in the introduction, this paper focuses on optimization problems of type \eqref{modprob_intro} 
where $\mathcal{R}$ is the risk measure given by \eqref{eq:risk}. 
In particular, we consider the conditional value-at-risk ($\mbox{CVaR}_{\beta}$) at confidence level 
$\beta, \; \beta\in (0,1)$ where $f$ in \eqref{eq:risk} is given by \eqref{eq:fcvar},
with $\mathcal{T}=\mathbb{R}$.	
CVaR$_{\beta}$ builds on the concept of value-at-risk (VaR) at level  $\beta \in (0,1)$, which is the $\beta$-quantile
of a given random variable.

More precisely, let $X$ be a random variable  and let  $\beta \in (0,1)$ be fixed. Then, VaR$_{\beta}[X]$ is given by
\[
\mbox{VaR}_{\beta}[X]:=\inf_{t\in\mathbb{R}} \, \{t:~\mathbb{P}[X\leq t]\geq \beta\},
\]
where  $\mathbb{P}[X\leq t]$ denotes the probability that the random variable $X$ is less than or
equal to $t$. 
VaR$_{\beta}$ is unfortunately not coherent because it violates sub-additivity / convexity. 
This is why CVaR$_{\beta}$ is preferred as a risk measure.

Even though now we know that coherent risk-measures can be written in the abstract form \eqref{eq:risk}, 
however, the 1D minimization formulation of CVaR$_{\beta}$ was first introduced  by Rockafellar and Uryasev 
in \cite{RU2000, RU2002}:
\[
\mbox{CVaR}_{\beta}[X]=\inf_{t\in\mathbb{R}} \left\{t + \frac{1}{1-\beta} \mathbb{E}[(X - t)_+]\right\}.
\]
 Moreover, if X is a continuous random variable, then
\[
 \mbox{CVaR}_{\beta}[X]=\mathbb{E}[X|X>\mbox{VaR}_{\beta}[X]],
\]
which shows that $\mbox{CVaR}_{\beta}[X]$ is the average of $\beta$-tail of the distribution of $X.$ Thus, $\mbox{CVaR}_{\beta}[X]$  emphasizes rare and low probability events, especially when $\beta  \rightarrow 1$.

\smallskip
\paragraph{\bf Model problem:} Our model problems are obtained by replacing $\mathcal{R}$ in \eqref{modprob_intro}
by $\mbox{CVaR}_{\beta}$, i.e., 
\be
\label{modprob}
\min_{u\in U_{ad}}\mbox{CVaR}_{\beta}[\mathcal{J}(y,u;\omega)] + \alpha\mathcal{P}(u) \;\; \mbox{subject\; to}\;\; c(y,u;\omega)=0, \quad \mbox{in } Z, \quad \mbox{a.a. } 
			\omega \in \Omega
\e
where again $U_{ad}$ denotes the set of admissible controls.
We use $(y, u)$ to denote the state and
control, respectively. In our setting $u$ is always deterministic.
The constraint equation $c(y,u;\omega)=0$ represents
 a differential equation with uncertain coefficients, $\mathcal{P}(u) $ is a deterministic cost 
 function, $\alpha$ is the regularization parameter, and $\mathcal{J}(y,u;\omega)$ 
 is a random variable cost function. 
 
 Here, we make the \emph{finite  dimensional noise assumption} on the equality constraint \cite{KS16}.
 We assume that $\omega$ can be sampled via a finite random vector $\xi : \Omega \rightarrow \Xi$ instead, where $ \Xi := \xi(\Omega)\subset \mathbb{R}^d$ with $d\in\mathbb{N}$.
 For example, coefficients, defining the constraint $c(y,u;\omega) = 0$, may be expressed by a Karhunen-Loeve (KL) approximation of an infinite-dimensional continuous random field, see \eqref{eq:kle} for an example.
This allows us to redefine the probability space to $(\Xi,\Sigma,\rho)$, where $\Sigma = \xi(\mathcal{A})$
is the $\sigma$-algebra of regions, and $\rho(\xi)$ is the continuous probability density function such that $\mathbb{E}[X] = \int_{\Xi} X(\xi)\rho(\xi) d\xi$.
The random variable $X(\xi)$ can be considered as a function of the random vector $\xi=(\xi^{(1)},\ldots,\xi^{(d)})$, belonging to the Hilbert space $\mathcal{F} = \{X(\xi): \|X\|<\infty\}$, equipped with the inner product $\langle X,Y\rangle = \int_{\Xi} X(\xi) Y(\xi) \rho(\xi) d\xi$ and the Euclidean norm $\|X\| = \sqrt{\langle X, X\rangle}$.
Since $\xi^{(1)},\ldots,\xi^{(d)}$ are independent random variables, we assume that each $\xi^{(k)}$ has a probability density function $\rho^{(k)}(\xi^{(k)})$, and that the space of functions $\mathcal{F}$ is isomorphic to a tensor product of spaces of univariate functions, $\mathcal{F} = \mathcal{F}^{(1)} \otimes \cdots \otimes \mathcal{F}^{(d)}$, where $\mathcal{F}^{(k)} = \{X^{{(k)}}(\xi^{(k)}): \|X^{(k)}\|<\infty\}$, $\|X^{(k)}\| = \sqrt{\langle X^{(k)}, X^{(k)} \rangle}$, $\langle X^{(k)}, Y^{(k)} \rangle = \int_{\mathbb{R}} X^{(k)}(\xi^{(k)}) Y^{(k)}(\xi^{(k)}) \rho^{(k)}(\xi^{(k)}) d\xi^{(k)}$, $k=1,\ldots,d$.

 Then, \eqref{modprob} reads
 \be
\label{modprob2}
\min_{u\in U_{ad}}\mbox{CVaR}_{\beta}[\mathcal{J}(y,u;\xi)] + \alpha\mathcal{P}(u) \;\; \mbox{subject\; to}\;\; c(y,u;\xi)=0,
 \quad \mbox{in } Z, \quad \mbox{a.a. } \xi \in \Xi .
\e
Here and in what follows, bracketed superscripts (e.g. $\xi^{(d)}$) denote a component, not a power or derivative.

To tackle nonsmoothness in CVaR$_\beta$, we employ a smoothing based approach from \cite{KS16}.
The smoothing approach is essentially aimed at approximating the positive part function $(\cdot)_+$ 
in $\mbox{CVaR}_{\beta}$ by a smooth function $g_{\varepsilon}: \mathbb{R} \rightarrow \mathbb{R}$, 
which depends on some $\varepsilon>0.$ Various examples of $g_{\varepsilon}$ are available in 
\cite[section~4.1.1]{KS16}. In particular, we consider the following 
$C^\infty$-smoothing  function
\begin{equation}
 \label{smoothf}
 g_{\varepsilon}(x) = \varepsilon \log(1+\exp(x/\varepsilon))
\end{equation}
where 
\begin{equation}
 \label{smoothf1}
 g'_{\varepsilon}(x) = \frac{1}{1+\exp(-x/\varepsilon)},
\qquad
g''_{\varepsilon}(x) = \frac{1}{\varepsilon} \left(\frac{1}{\exp(x/(2\varepsilon)) + \exp(-x/(2\varepsilon))}\right)^2.
\end{equation}

Thus, the optimization problem for smooth CVaR$_\beta^\varepsilon$ is given by
\be
\label{modprob_sm}
\begin{cases}
\min\limits_{(u,t)\in U_{ad} \times \mathbb{R}}
	 \mathcal{R}_{t,\beta}^\varepsilon[\mathcal{J}(y,u;\xi)]  + \alpha\mathcal{P}(u)  \\
\mbox{subject\; to} \\
c(y,u;\xi)=0, \quad \mbox{in } Z, \quad \mbox{a.a. } \xi \in \Xi , 
			\end{cases}
\e
where 
\begin{equation}\label{eq:cvar_smooth}
	\mathcal{R}^{\varepsilon}_{t,\beta}[\mathcal{J}(y,u;\xi)]
	:= t + \frac{1}{1-\beta} \mathbb{E}[g_{\varepsilon}(\mathcal{J}(y,u;\xi) - t)] .
\end{equation}
For convergence analysis of \eqref{modprob_sm} to \eqref{modprob2}, we refer 
to \cite{KS16}.

\subsection{Cartesian function space}
The dimension $d$ of the random vector $\xi$ can be arbitrarily high, e.g. tens of model tuning parameters or KL coefficients.
In this case expectations as in \eqref{eq:cvar_smooth} become high-dimensional integrals.
Instead of a direct Monte Carlo average (which may converge too slowly),
we can introduce a high-order quadrature rule (e.g. Gauss-Legendre) with $n_{\xi} \in \mathbb{N}$ points in each of the components $\xi^{(1)},\ldots,\xi^{(d)}$ independently.
However, the exponential total number of quadrature points in all variables $n_{\xi}^d$ becomes intractable even for moderate dimensions.

\subsection{Tensor Train decomposition}
We circumvent this ``curse of dimensionality'' problem by approximating all functions depending on $\xi$ by a Tensor Train (TT) decomposition~\cite{osel-tt-2011},
which admits efficient integration and differentiation.
\begin{definition}
 A square-integrable function $f(\xi)$ is said to be approximated by a \emph{(functional) TT decomposition} \cite{Marzouk-stt-2016,Gorodetsky-ctt-2019} with a relative approximation error $\epsilon$ if there exist univariate functions $F^{(k)}(\cdot):\xi^{(k)} \in \mathbb{R} \rightarrow \mathbb{R}^{r_{k-1}\times r_k}$, $k=1,\ldots,d$, such that
 \begin{equation}\label{eq:ftt}
  \tilde f(\xi) := \sum_{s_0,\ldots,s_d=1}^{r_0,\ldots,r_d} F^{{(1)}}_{s_0,s_1}(\xi^{(1)}) F^{{(2)}}_{s_1,s_2}(\xi^{(2)}) \cdots F^{{(d)}}_{s_{d-1},s_d}(\xi^{(d)}),
 \end{equation}
 where the subscripts $s_{k-1},s_k$ denote elements of a matrix,
 and $\|f - \tilde f\| = \epsilon \|f\|$.
 The factors $F^{(k)}$ are called \emph{TT cores}, and the ranges of summation indices $r_0,\ldots,r_d \in \mathbb{N}$ are called \emph{TT ranks}.
\end{definition}

Without loss of generality we can let $r_0=r_d=1$, but the other TT ranks $r_1, \ldots,r_{d-1}$  can vary depending on the approximation error.
One example is a bi-variate truncated Fourier series $\tilde f(\xi^{(1)},\xi^{(2)}) = \sum_{s=-r}^{r} f_s(\xi^{(1)}) \exp(\mathrm{i}s \xi^{(2)})$.

From \eqref{eq:ftt}, we notice that the expectation of $\tilde f$ factorizes into univariate integrations,
$$
\mathbb{E} [\tilde f] = \sum_{s_0,\ldots,s_d=1}^{r_0,\ldots,r_d} \left(\int F^{(1)}_{s_0,s_1}(\xi^{(1)}) \rho^{(1)}(\xi^{(1)}) d\xi^{(1)} \right) \cdots \left(\int F^{(d)}_{s_{d-1},s_d}(\xi^{(d)}) \rho^{(d)}(\xi^{(d)}) d\xi^{(d)} \right).
$$

For practical computations with \eqref{eq:ftt} we introduce univariate bases $\{\ell_i(\xi^{(k)})\}_{i=1}^{n_{\xi}}$,
and the multivariate basis constructed from a tensor product,
$$
L_{i_1,\ldots,i_d}(\xi):= \ell_{i_1}(\xi^{(1)}) \cdots \ell_{i_d}(\xi^{(d)}).
$$
Now we can collect the expansion coefficients of $\tilde f$ into a tensor $\mathbf{F}\in\mathbb{R}^{n_{\xi} \times \cdots \times n_{\xi}}$,
\begin{equation}\label{eq:basis}
 \tilde f(\xi) = \sum_{i_1,\ldots,i_d=1}^{n_{\xi}} \mathbf{F}(i_1,\ldots,i_d) L_{i_1,\ldots,i_d}(\xi).
\end{equation}
Similarly, TT cores in \eqref{eq:ftt} can be written using three-dimensional tensors $\mathbf{F}^{(k)}\in\mathbb{R}^{r_{k-1} \times n_{\xi} \times r_k}$,
\begin{equation}
 F^{(k)}_{s_{k-1},s_k}(\xi^{(k)}) = \sum_{i=1}^{n_{\xi}} \mathbf{F}^{(k)}(s_{k-1},i,s_k) \ell_i(\xi^{(k)}), \quad k=1,\ldots,d.
\end{equation}
The original (discrete) TT decomposition \cite{osel-tt-2011} was introduced for tensors,
\begin{equation}
 \mathbf{F}(i_1,\ldots,i_d) = \sum_{s_0,\ldots,s_d=1}^{r_0,\ldots,r_d} \mathbf{F}^{(1)}(s_0,i_1,s_1) \cdots \mathbf{F}^{(d)}(s_{d-1},i_d,s_d).
\end{equation}
Note that $\mathbf{F}$ contains $n_{\xi}^d$ elements, whereas storing $\mathbf{F}^{(1)},\ldots,\mathbf{F}^{(d)}$ needs only $\sum_k r_{k-1} n_\xi r_k$ elements.
For brevity we can define the maximal TT rank $r:=\max_k r_k$,
which gives us a linear storage complexity of the TT decomposition, $\mathcal{O}(dn_{\xi} r^2)$.

If $\{\ell_i(\xi^{(k)})\}_{i=1}^{n_{\xi}}$ is a Lagrange polynomial basis, defined by the quadrature points $\Xi^{(k)}:=\{\xi^{(k)}_{i}\}$ and weights $\{w_i\}$ such that $\ell_i(\xi^{(k)}_j) = \delta_{i,j}$, we obtain
\begin{align}\label{eq:quad}
 \mathbb{E} [\tilde f] & = \sum_{i_1,\ldots,i_d} \tilde f(\xi_{i_1}^{(1)},\ldots,\xi_{i_d}^{(d)}) w_{i_1}\cdots w_{i_d} =  \sum_{i_1,\ldots,i_d} \mathbf{F}(i_1,\ldots,i_d) w_{i_1}\cdots w_{i_d}\\
 & = \sum_{s_0,\ldots,s_d=1}^{r_0,\ldots,r_d} \left(\sum_{i_1} \mathbf{F}^{(1)}(s_0,i_1,s_1) w_{i_1} \right) \cdots \left(\sum_{i_d} \mathbf{F}^{(d)}(s_{d-1},i_d,s_d) w_{i_d} \right). \nonumber
\end{align}
The summation over $s_0,\ldots,s_d$ in the right hand side can be computed recursively by multiplying only two tensors at a time.
Assuming that a partial result $\mathbf{R}_k\in\mathbb{R}^{r_0 \times r_{k}}$ is given, we can compute
\begin{equation}\label{eq:int-rec}
\mathbf{R}_{k+1}(s_0,s_{k+1}) = \sum_{s_k=1}^{r_k} \mathbf{R}_k(s_0,s_k) \left(\sum_{i_{k+1}} \mathbf{F}^{(k+1)}(s_{k},i_{k+1},s_{k+1}) w_{i_{k+1}} \right)
\end{equation}
as a matrix product with a $\mathcal{O}(n_\xi r^2)$ complexity.
Starting with $\mathbf{R}_0=1$ and finishing with $\mathbf{R}_d = \mathbb{E}[\tilde f]$, we complete
the entire integration in $\mathcal{O}(dn_{\xi} r^2)$ operations.

Such a recursive sweep over TT cores is paramount to computing TT approximations of arbitrary functions, or to solution of operator equations in the TT format.
For example, the \emph{TT-Cross} method (see Appendix~\ref{sec:cross}) requires $\mathcal{O}(dn_{\xi}r^2)$ samples from a function $f(\xi)$ and $\mathcal{O}(dn_{\xi}r^3)$ further floating point operations to compute a TT approximation $\tilde f(\xi) \approx f(\xi)$.
Similarly, linear algebra on functions can be recast to linear algebra on their TT cores with a linear complexity in the dimension (see Appendix~\ref{sec:ttalg}).

\section{Reduced Space Formulation}\label{sec:reduced}

This paper considers two approaches to tackle \eqref{modprob_sm}.
In this section, we first introduce TTRISK for the so-called reduced form of \eqref{modprob_sm}.
To ensure that CVaR$_\beta^\varepsilon$ is a statistically unbiased estimator of $\mbox{CVaR}_{\beta}$, 
we introduce a control variate correction in subsection~\ref{s:control_variate}.

\subsection{Smoothed CVaR with TT approximations}
Assume that $c(y,u;\xi) = 0$ is uniquely solvable,
i.e., for each $u \in U_{ad}$ there exists a unique solution mapping 
$y(u;\cdot) : \Xi \rightarrow Y$ for $\mathbb{P}$ a.a $\xi \in \Xi$. The 
resulting optimization problem \eqref{modprob_sm} only depends on $u$ and 
is given by 
\be
\label{modprob_redsm}
\min\limits_{(u,t)\in U_{ad} \times \mathbb{R}} \left\{ \mathfrak{J}(u,t)
	:= \mathcal{R}^{\varepsilon}_{t,\beta}[j(u;\xi)]  + \alpha\mathcal{P}(u) \right\} ,
\e
where 
\begin{equation}\label{eq:jred}
	j(u;\xi) := \mathcal{J}(y(u;\xi),u;\xi) .
\end{equation}

The exact expectation in $\mbox{CVaR}^{\varepsilon}_\beta$
can be approximated by a quadrature similarly to \eqref{eq:quad}.
We denote the total number of quadrature points $N$ (which is formally $N=n_{\xi}^d$ in the Cartesian formulation).
However, we need to tackle the curse of dimensionality using the TT decomposition.

\begin{definition}
 The approximate expectation $\mathbb{E}_N [f]$ of a function $f(\xi)$ is defined as $\mathbb{E}_N [f] := \mathbb{E} [\tilde f]$, where $\tilde{f}(\xi)$ is a TT approximation \eqref{eq:ftt} to $f(\xi)$, computed using the TT-Cross algorithm as described in Appendix~\ref{sec:cross}, and the integration of $\tilde f(\xi)$ is carried out as shown in \eqref{eq:int-rec}.
\end{definition}

This leads to the following approximation of \eqref{modprob_redsm}:
\be
\label{modprob_red_disc}
\min\limits_{(u,t) \in \mathcal{U}_{ad} \times \mathbb{R}}  \left\{ \mathfrak{J}_N(u,t)
 := \mathcal{R}_{t,\beta,N}^\varepsilon[j(u;\xi)]  + \alpha\mathcal{P}(u) \right\} ,
\e
where 
\[
	\mathcal{R}_{t,\beta,N}^\varepsilon[j(u;\xi)]
	:= t + \frac{1}{(1-\beta)} \mathbb{E}_N[g_{\varepsilon}(j(u;\xi) - t)].
\]
However, to optimize the entire cost $\mathfrak{J}_N$
we need
to calculate the first and second order derivatives. We readily obtain that the first order derivatives
are
\begin{equation}\label{eq:Ju}
\begin{aligned}
	\nabla_u \mathfrak{J}_N(u,t) 
	 &= (1-\beta)^{-1} \mathbb{E}_N[g'_\varepsilon(j(u;\xi)-t) \nabla_u j(u;\xi)] + \alpha \nabla_u \mathcal{P}(u) \\
	\nabla_t \mathfrak{J}_N(u,t) 
	 &= 1-(1-\beta)^{-1} \mathbb{E}_N[g'_\varepsilon(j(u;\xi)-t) ] 
\end{aligned}
\end{equation}
and the second order derivatives are
\begin{align}\nonumber
	\nabla_{uu} \mathfrak{J}_N(u,t)  
	 &= (1-\beta)^{-1} \mathbb{E}_N \left[g''_\varepsilon(j(u;\xi)-t) \nabla_u j(u;\xi) \nabla_u j(u;\xi)^* 
                + g'_\varepsilon(j(u;\xi)-t) \nabla_{uu} j(u;\xi) \right] \\ \label{eq:Juu}
	 &\quad  + \alpha \nabla_{uu} \mathcal{P}(u) , \\\label{eq:Jut}
	\nabla_{ut} \mathfrak{J}_N(u,t)  
	 &= - (1-\beta)^{-1} \mathbb{E}_N[g''_\varepsilon(j(u;\xi)-t) \nabla_u j(u;\xi)] , \\
	\nabla_{tu} \mathfrak{J}_N(u,t)    
	 &= - (1-\beta)^{-1} \mathbb{E}_N[g''_\varepsilon(j(u;\xi)-t) \nabla_u j(u;\xi)^*] , \\
	\nabla_{tt} \mathfrak{J}_N(u,t)     
	 &= (1-\beta)^{-1} \mathbb{E}_N[g''_\varepsilon(j(u;\xi)-t) ]. \label{eq:Jtt}
\end{align}
Observe that the second derivatives of $\mathfrak{J}_N$ computed above depend on  $g''_{\varepsilon}.$ 
From \eqref{smoothf1}, we see that $g''_{\varepsilon}$ decays rapidly away from the interval 
$(-\varepsilon, \varepsilon)$.
Consequently, the Hessian
\begin{equation}\label{eq:H}
 H =
 \begin{bmatrix}
 \nabla_{uu} \mathfrak{J}_N & \nabla_{ut} \mathfrak{J}_N\\
 \nabla_{tu} \mathfrak{J}_N & \nabla_{tt} \mathfrak{J}_N
 \end{bmatrix}
\end{equation}
can degenerate if $g''_{\varepsilon}(j(u;\xi_i) - t)=0$ for all $i=1,\ldots, N,$  since all but
the $(1,1)$ block of $H$ are zeros \cite{MaterThesisMae}.
(In fact, the Hessian may become ill-conditioned also if $g''_{\varepsilon}(j(u;\xi_i) - t)$ is close to zero.)
To circumvent this problem, we adopt a technique similar to that used in augmented Lagrangian methods to dynamically update augmentation parameters within the optimization problem \cite{fgg,MaterThesisMae}.
We start with $t_0 = \varepsilon_0 = \mathbb{E}_N [j(u_0;\xi)]$.
In many cases there will be a lot of points $j(u;\xi_i)-t$ in a significant support of $g''_{\varepsilon}$.
During the course of Newton iterations, we decrease $\varepsilon$ geometrically,
\begin{equation}\label{eq:mu_iter}
 \varepsilon_{m+1} = \mu_{\varepsilon} \varepsilon_{m},
\end{equation}
where $0<\mu_{\varepsilon}<1$ is a tuning factor (e.g. $\mu_{\varepsilon}=1/2$).
However, the next iterate $t_{m+1}$ may end up far on the tail of $g''_{\varepsilon_m}$ again.
To prevent this from happening, we perform a line search, where in addition to
a non-increasing residual condition\footnote{See Remark~\ref{rem:armijo} for further explanation of this condition.}
\begin{equation}\label{eq:Arm1}
 \left\|\begin{bmatrix}\grad{\mathfrak{J}_N}{u}(u_{m+1},t_{m+1}) \\ \grad{\mathfrak{J}_N}{t}(u_{m+1},t_{m+1})\end{bmatrix}\right\|  \le \left\|\begin{bmatrix}\grad{\mathfrak{J}_N}{u}(u_m,t_m) \\ \grad{\mathfrak{J}_N}{t}(u_m,t_m)\end{bmatrix}\right\|,
\end{equation}
where $u_{m+1} = u_m + h \delta_{um}$, $t_{m+1} = t_m + h \delta_{tm}$ with a step size $h>0$ and Newton directions $\delta_{um},\delta_{tm}$, we require that
\begin{equation}\label{eq:Arm2}
\mathbb{E}_N [\exp(-|j(u_{m+1};\xi)-t_{m+1}|/\varepsilon_m)] > \theta
\end{equation}
for some $0<\theta<1$.
In the numerical experiments we have found the iterations to be robust and insensitive for $\theta$ between $10^{-2}$ and $10^{-1}$.
\begin{proposition}
This ensures that
$$
\nabla_{tt} \mathfrak{J}_N(u_{m+1},t_{m+1}) > \frac{\theta^{1/\mu_{\varepsilon}}}{4\varepsilon_{m+1} (1-\beta)}
$$
stays away from zero by a fixed fraction of the maximum of $g''_{\varepsilon_{m+1}}(t)$, which is equal to $1/(4\varepsilon_{m+1})$.
\end{proposition}
\begin{proof}
Firstly, note from \eqref{smoothf1} that
\begin{align}\label{eq:gppmax}
g''_{\varepsilon}(t) & \ge \frac{1}{\varepsilon} \left(\frac{1}{2\exp(|t|/(2\varepsilon))}\right)^2 = \frac{1}{4\varepsilon} \exp\left(-\frac{|t|}{\varepsilon}\right). \tag{*}
\end{align}
Now,
\begin{align*}
\nabla_{tt} \mathfrak{J}_N(u_{m+1},t_{m+1}) & = (1-\beta)^{-1} \mathbb{E}_N \left[g''_{\varepsilon_{m+1}}(j(u_{m+1};\xi)-t_{m+1})\right]  & \mbox{(from Eq. \eqref{eq:Jtt})}\\
& \ge  \frac{1}{4\varepsilon_{m+1} (1-\beta)} \mathbb{E}_N \left[ \exp\left(-\frac{|j(u_{m+1};\xi)-t_{m+1}|}{\mu_{\varepsilon}\varepsilon_m}\right)\right] & \mbox{(by linearity of $\mathbb{E}_N$, \eqref{eq:gppmax} and \eqref{eq:mu_iter})} \\
& =  \frac{1}{4\varepsilon_{m+1} (1-\beta)} \mathbb{E}_N \left[ \exp\left(-\frac{|j(u_{m+1};\xi)-t_{m+1}|}{\varepsilon_m}\right)^{1/\mu_\varepsilon}\right] \\
& \ge \frac{1}{4\varepsilon_{m+1} (1-\beta)} \left(\mathbb{E}_N \left[ \exp\left(-\frac{|j(u_{m+1};\xi)-t_{m+1}|}{\varepsilon_m}\right)\right] \right)^{1/\mu_\varepsilon} & \mbox{(by Jensen's inequality since $\frac{1}{\mu_{\varepsilon}}>1$)} \\
& > \frac{1}{4\varepsilon_{m+1} (1-\beta)} \theta^{1/\mu_{\varepsilon}}. & \mbox{(by assumption \eqref{eq:Arm2})}
\end{align*}
The proof is complete.
\end{proof}
If the dimension of discretized $u$ is small, we can compute the TT-Cross approximation of $\grad{j}{u}(u_m,\xi),\grad{j}{uu}(u_m,\xi),g''_{\varepsilon_m}(j(u_m;\xi)-t_m)$ and $g''_{\varepsilon_m}(j(u_m;\xi)-t_m)$ directly,
since evaluations of $j$, $\grad{j}{u}$ and $\grad{j}{uu}$ are available explicitly from \eqref{eq:jred} under the unique solution mapping
$y(u;\cdot)$,
compute the expectations in \eqref{eq:Ju} and \eqref{eq:Juu},
and solve
\begin{equation}\label{eq:NewtonSys}
 H \begin{bmatrix}\delta_{um} \\ \delta_{tm} \end{bmatrix} = -\begin{bmatrix}\grad{\mathfrak{J}_N}{u}(u_m,t_m) \\ \grad{\mathfrak{J}_N}{t}(u_m,t_m)\end{bmatrix}
\end{equation}
for the Newton directions.
However, if $u$ is large, $\grad{j}{uu}(u_m,\xi)$ is large and dense, and its TT decomposition becomes too expensive.
In fact, even multiplying $\grad{\mathfrak{J}_N(u,t)}{uu}$ by a vector in an iterative solver would require recomputation of the TT decomposition in each iteration.
To avoid this problem, we propose to replace the expectation in \eqref{eq:Juu} by sampling at a fixed point $\bar\xi$.
This gives
\begin{align}
\grad{\widetilde{\mathfrak{J}_N}}{uu}(u,t)
	 & = (1-\beta)^{-1} \left[\mathbb{E}_N[g''_\varepsilon(j(u;\xi)-t)] \nabla_u j(u;\bar\xi) \nabla_u j(u;\bar\xi)^* \right] \nonumber \\
               & \quad + (1-\beta)^{-1} \mathbb{E}_N[g'_\varepsilon(j(u;\xi)-t)] \nabla_{uu} j(u;\bar\xi)  + \alpha \nabla_{uu} \mathcal{P}(u), \label{eq:JuuFixed}
\end{align}
and consequently a \emph{fixed-point Hessian}
\begin{equation}\label{eq:HFixed}
\widetilde{H} =
 \begin{bmatrix}
 \nabla_{uu} \widetilde{\mathfrak{J}_N} & \nabla_{ut} \mathfrak{J}_N\\
 \nabla_{tu} \mathfrak{J}_N & \nabla_{tt} \mathfrak{J}_N
 \end{bmatrix} . 
\end{equation}
The choice of $\bar\xi$ is motivated by the mean value theorem.
We can treat $\mathbb{E}_N[g'_\varepsilon(j(u;\xi)-t) \nabla_{uu} j(u;\xi)]$ as an expectation of $\nabla_{uu} j(u;\xi)$ alone over a probability density function $\bar\rho(\xi) = (1/C) g'_\varepsilon(j(u;\xi)-t) \rho(\xi)$, where $C := \mathbb{E}_N[g'_\varepsilon(j(u;\xi)-t)]$ is the normalizing constant. In turn, for a linear $\nabla_{uu} j(u;\xi)$ it holds $\mathbb{E}_{\bar\rho}[\nabla_{uu} j(u;\xi)] = \nabla_{uu} j(u;\mathbb{E}_{\bar\rho}[\xi])$,
so we can take the right hand side as an approximation also in a general case.
This gives
\begin{equation}\label{eq:xiFixed}
\bar\xi = \mathbb{E}_{\bar\rho} [\xi] = \frac{\mathbb{E}_N[g'_\varepsilon(j(u;\xi)-t) \xi]}{\mathbb{E}_N[g'_\varepsilon(j(u;\xi)-t)]}.
\end{equation}
We focus on $\nabla_{uu} j$ (and hence on $g'_\varepsilon$) since $\nabla_{uu} j$ is usually the dominant part of $\nabla_{uu} \mathfrak{J}_N$.
Note that the action $\nabla_{uu} j(u;\bar\xi) \cdot \delta_u$ can usually be applied efficiently, since this requires the solution of one forward and one adjoint deterministic problems at fixed $\xi=\bar\xi$.
Similarly, $\nabla_{uu} \mathcal{P}(u)$ is a sparse, and $\nabla_u j(u;\bar\xi) \nabla_u j(u;\bar\xi)^*$ is a rank-1 matrix after discretization.
This allows us to solve the Newton system \eqref{eq:NewtonSys} with $\widetilde{H}$ instead of $H$ iteratively with fast matrix-vector products.

Lastly, if $U_{ad}$ is constrained, we can add the projection of the Newton direction onto $U_{ad}$.
To reduce the number of projections in the step selection stage, we write the method in a Frank-Wolfe's fashion \cite{bertsekas1999nonlinear}, that is, we project the search direction,  
$
\hat\delta_{um} := \mathrm{Proj}_{U_{ad}}(u_m+\delta_{um}) - u_m,
$
where $\mathrm{Proj}_{U_{ad}}(\cdot)$ is the orthogonal projection onto $U_{ad}$,
followed by the usual line search in $u_{m+1} = u_m + h \hat\delta_{um}$.
The entire procedure is summarized in Algorithm~\ref{alg:ttrisk}.


 \begin{algorithm}[htb]
\caption{TTRISK (Reduced Case)}
\label{alg:ttrisk}
 \begin{algorithmic}[1]
  \Require $\alpha$, $\beta$, line search parameter $\theta$, number of iterations $I_{\max}$, smoothness reduction factor $\mu_{\varepsilon}$, stopping tolerance, procedure to compute $j(u,\xi)$.
  \State Set $m=0$, $u_0=0$, $t_0 = \varepsilon_0 = \mathbb{E}_N [j(u_0;\xi)]$.
  \While{$m\le I_{\max}$ and $|t_m-t_{m-1}| > \mbox{tol} \cdot |t_m|$ or $\|u_m-u_{m-1}\| > \mbox{tol} \cdot \|u_m\|$ or $m=0$}
   \State Approximate $j(u_m;\xi),\grad{j(u_m;\xi)}{u},g'_{\varepsilon_m}(j(u_m;\xi)-t_m),g''_{\varepsilon_m}(j(u_m;\xi)-t_m)$ by TT-Cross.
   \State Compute expectations in \eqref{eq:Ju}, \eqref{eq:Jut}--\eqref{eq:Jtt}, \eqref{eq:JuuFixed} and \eqref{eq:xiFixed} using \eqref{eq:quad},\eqref{eq:int-rec}.
   \State Solve \eqref{eq:NewtonSys} with $H$ or $\widetilde{H}$ using Conjugate Gradients method.
   \State Project the increment $\hat \delta_{um} := \mathrm{Proj}_{U_{ad}}(u_m+\delta_{um}) - u_m$. \Comment{{for control constraints}}
   \State Find $0<h\le 1$ such that \eqref{eq:Arm1} and \eqref{eq:Arm2} hold.
   \State Update $u_{m+1} = u_m + h \hat\delta_{um}$, $t_{m+1} = t_m + h \delta_{tm}$, $\varepsilon_{m+1} = \mu_{\varepsilon} \varepsilon_m$.
   \State Set $m=m+1$
  \EndWhile\\
  \Return $u_m\approx u$, $t_m\approx t$, $\mathcal{R}_{t,\beta,N}^{\varepsilon_m}[j(u_m;\xi)] \approx \mbox{CVaR}_{\beta}$.
 \end{algorithmic}
\end{algorithm}


\subsection{Smoothed CVaR as control variate for Monte Carlo}
\label{s:control_variate}

Assuming that $\tilde g_{\varepsilon}$ and $\tilde j$ are TT approximations to $g_{\varepsilon}$ and $j$, respectively,
we can define the approximate expectation of the exact function as the exact expectation of the approximate function, since the approximate function is a polynomial:
$$
\mathbb{E}_N[g_{\varepsilon}(j(u;\xi) - t)] = \mathbb{E}[\tilde g_{\varepsilon}(\tilde j(u;\xi) - t)].
$$
Using add-and-subtract trick, we can write the exact risk-measure functional as follows,
$$
\mathcal{R}_t = t + \frac{1}{1-\beta} \mathbb{E}[\tilde g_{\varepsilon}(\tilde j(u;\xi) - t)] + \frac{1}{1-\beta} \left[\mathbb{E}(j(u;\xi) - t)_+ - \mathbb{E}[\tilde g_{\varepsilon}(\tilde j(u;\xi) - t)]\right].
$$
However, for the last term we can use a different quadrature for computing the expectation, such as Monte Carlo with $M$ samples.
This defines a \emph{corrected} smoothed functional:
\begin{align}
\mathcal{R}_{t,\beta,N}^{\varepsilon,M}[j(u;\xi)]
	& := t + \frac{1}{(1-\beta)} \mathbb{E}_N[\tilde g_{\varepsilon}(\tilde j(u;\xi) - t)] \nonumber \\
	& + \frac{1}{(1-\beta)}\frac{1}{M} \sum_{\ell=1}^{M} \left[(j(u;\xi_{\ell}) - t)_+ - \tilde g_{\varepsilon}(\tilde j(u;\xi_{\ell}) - t)\right], \label{eq:MCred}
\end{align}
where $\xi_{\ell}$ are i.i.d. samples from $\rho(\xi)$.
The benefit of this scheme stems from the fact that if the approximation is accurate, i.e. $\mathrm{var}[(j(u;\xi) - t)_+ - \tilde g_{\varepsilon}(\tilde j(u;\xi) - t)] \le \delta^2 \mathrm{var}[(j(u;\xi) - t)_+]$ is small (where $\delta$ is the total error from the TT-Cross (Algorithm~\ref{alg:tt-cross}) and smoothing),
by the law of large numbers
the variance of this estimator reads
$$
\mathrm{var}[\mathcal{R}_{t,\beta,N}^{\varepsilon,M}] \le \delta^2 \frac{\mathrm{var}[(j(u;\xi) - t)_+]}{M} \ll \frac{\mathrm{var}[(j(u;\xi) - t)_+]}{M},
$$
where the latter term is the variance of the straightforward Monte Carlo approximation of $\mathcal{R}_{t}$.
Consequently, one needs a much smaller $M$ to achieve the same variance (i.e. Mean Square Error) threshold.
One can say that $\tilde g_{\varepsilon}$ is used as a \emph{control variate} for variance reduction of Monte Carlo \cite{robert-MC-book-2004,nobile-mlmcwithcv-2015}, or vice versa, that \eqref{eq:MCred} is the \emph{2nd level correction} \cite{giles-mlmc-2015,Scheichl-mlqmc-lognorm-2017} to the surrogate model $\tilde g_{\varepsilon}$.
This makes $\mathcal{R}_{t,\beta,N}^{\varepsilon,M}$ a statistically unbiased estimator of $\mathcal{R}_{t}$.

Similarly we can correct the cost function and its gradient \eqref{eq:Ju}:
\begin{align}
	\nabla_u \mathfrak{J}_N^M(u,t)
	 &= (1-\beta)^{-1} \mathbb{E}[\tilde g'_\varepsilon(\tilde j(u;\xi)-t) \nabla_u \tilde j(u;\xi)] + \alpha \nabla_u \mathcal{P}(u) \nonumber \\
	 & + \frac{1}{(1-\beta)M}\sum_{\ell=1}^M \left[\theta(j(u;\xi_\ell)-t) \nabla_u j(u;\xi_\ell) - \tilde g'_\varepsilon(\tilde j(u;\xi_\ell)-t) \nabla_u \tilde j(u;\xi_\ell) \right], \label{eq:MCgradu}\\
	\nabla_t \mathfrak{J}_N^M(u,t)
	 &= 1-(1-\beta)^{-1} \mathbb{E}[\tilde g'_\varepsilon(\tilde j(u;\xi)-t) ] \nonumber \\
	 & - \frac{1}{(1-\beta)M} \sum_{\ell=1}^{M}\left[\theta(j(u;\xi_\ell)-t) - \tilde g'_\varepsilon(\tilde j(u;\xi_\ell)-t)\right], \label{eq:MCgradt}
\end{align} 
where $\theta(t)=1$ if $t\ge 0$, and $0$, otherwise.
These gradients can be plugged into Line~4 of Algorithm~\ref{alg:ttrisk} instead of \eqref{eq:Ju}.
Since the variance of the correction is expected to be small,
we omit it in the Hessian, turning Alg.~\ref{alg:ttrisk} into a Gauss-Newton method.

\section{Lagrangian CVaR formulation}\label{sec:lagrangian}

In this section, we first focus on the full-space formulation. This is followed by a Gauss-Newton system setup for
the problem and a preconditioning strategy for this system.

\subsection{Semi-discrete formulation}
As in Section~\ref{sec:reduced}, we let $\{\xi_i\}_{i=1}^N$ be all nodes, and $\{w_i\}_{i=1}^{N}$ be all weights in the quadrature \eqref{eq:quad}.
Moreover, we assume a Lagrangian basis expansion \eqref{eq:basis}, and for brevity we let $L_i(\xi):=L_{i_1,\ldots,i_d}(\xi)$ and $f_i := \tilde f(\xi_i) = \mathbf{F}(i_1,\ldots,i_d)$.
Applying this formalism to $y$ instead of $f$, and
assuming that the constraint $c(y,u,\xi)=0$ holds pointwise in $\xi$,
we obtain semi-discrete equations
\[
 c(y_i, u, \xi_i)=0, \qquad i=1,\ldots,N.
\]
Likewise, we can introduce the smoothed CVaR with the Monte Carlo correction (cf.~\eqref{eq:MCred}) using the quadrature
\begin{align*}
 \mathcal{R}^{\varepsilon,M}_{t,\beta,N} & = t + (1-\beta)^{-1} \sum_{i=1}^{N} w_i g_{\varepsilon}(\mathcal{J}(y_i,u,\xi_i) - t) \\
&\quad + \frac{1}{(1-\beta)M}\sum_{\ell=1}^{M} \left[(\mathcal{J}(\tilde y(\xi_\ell),u,\xi_\ell)-t)_+ - g_{\varepsilon}(\mathcal{J}(\tilde y(\xi_\ell),u,\xi_\ell) - t)\right].
\end{align*}
Let $\vec{y}\in Y^{N}$ and $\vec{p}\in (Y^*)^{N}$ be functions $y_i$ and $p_i$ stacked together.
We introduce the Lagrangian
\begin{align*}
\mathcal{L}(\vec{y},u,\vec{p},t) & = \mathcal{R}^{\varepsilon,M}_{t,\beta,N} + \alpha \mathcal{P}(u) + \sum_{i=1}^{N} w_i \langle p_i, c(y_i,u,\xi_i) \rangle.
\end{align*}
In the differentiation, we will distinguish the components corresponding to different coefficients.
This gives, for all $j=1,\ldots,N$,
\begin{align}
\grad{\mathcal{L}}{y_j} & = \frac{1}{1-\beta} \left[w_j g_{\varepsilon}' (\mathcal{J}(y_j,u,\xi_j) - t) \grad{\mathcal{J}}{y}(y_j,u,\xi_j) + \frac{1}{M} \sum_{\ell=1}^{M} E_j'(\xi_\ell)\right]  + w_j (\grad{c}{y})^* p_j, \nonumber \\
\grad{\mathcal{L}}{u} & = \alpha \grad{\mathcal{P}}{u}(u) + \sum_{i=1}^{N} w_i (\grad{c}{u}(y_i,u,\xi_i))^* p_i , \label{eq:dLdu}\\
\grad{\mathcal{L}}{p_j} & = w_j c(y_j,u,\xi_j), \nonumber \\
\grad{\mathcal{L}}{t} & = 1 - (1-\beta)^{-1} \left[\sum_{i=1}^{N} w_i g_{\varepsilon}'(\mathcal{J}(y_i,u,\xi_i) - t) + \frac{1}{M}\sum_{\ell=1}^{M} e'(\xi_\ell)\right], \nonumber
\end{align}
where we have defined error correction shortcuts
\begin{align}
e'(\xi) & = \left[\theta\left(\mathcal{J}\left(\tilde y(\xi),u,\xi\right)-t\right)  - g_{\varepsilon}'\left(\mathcal{J}\left(\tilde y(\xi),u,\xi\right)-t\right)\right], \nonumber\\
E_j'(\xi) & = e'(\xi) L_j(\xi) \grad{\mathcal{J}}{y} \left(\tilde y(\xi),u,\xi\right) ,
\label{eq:corrderiv}
\end{align}
with $\theta(t)=1$ for $t\ge 0$ and $0$ otherwise.

The second derivatives of \eqref{eq:corrderiv} are difficult both notationally and computationally, since arbitrary points $\xi_\ell$, leading to nonzero Lagrangian polynomial values $L_j(\xi_\ell)$, produce dense matrices.
However, if we assume that the correction is small in magnitude, we can follow the arguments of Gauss-Newton methods and remove the correction derivatives in the Hessian, as well as second derivatives of $c(y,u,\xi)$ and $\mathcal{J}(y,u,\xi)$.
This way we obtain (where $\delta_{j,k}$ denotes the Kronecker delta, and $[f(y,u,\xi)]_j$ evaluates $f(y_j,u,\xi_j)$)
\begin{align*}
 \grad{\mathcal{L}}{y_j, y_k} & \approx (1-\beta)^{-1} \delta_{j,k} w_j \left[g_{\varepsilon}'' \grad{\mathcal{J}}{y} (\grad{\mathcal{J}}{y})^* + g_{\varepsilon}' \grad{\mathcal{J}}{yy} \right]_{j}, & j,k=1,\ldots,N \\
 \grad{\mathcal{L}}{y_j, p_k} & \approx \delta_{j,k} w_j (\grad{c(y_j,u,\xi_j)}{y})^*, & j,k=1,\ldots,N \\
 \grad{\mathcal{L}}{y_j, t} & \approx - (1-\beta)^{-1} w_j g_{\varepsilon}''(\mathcal{J}(y_j,u,\xi_j)-t) \grad{\mathcal{J}(y_j,u,\xi_j)}{y}, & j=1,\ldots,N \\[2ex]
 \grad{\mathcal{L}}{u, p_k} & \approx (\grad{c(y_k,u,\xi_k)}{u})^* w_k, & k=1,\ldots,N \\
 \grad{\mathcal{L}}{p_j, y_k} & \approx \delta_{j,k} w_j \grad{c(y_j,u,\xi_j)}{y}, & j,k=1,\ldots,N \\
 \grad{\mathcal{L}}{p_j, u} & \approx w_j \grad{c(y_j,u,\xi_j)}{u}, & j=1,\ldots,N\\
 \grad{\mathcal{L}}{t,t} & \approx (1-\beta)^{-1} \sum_{i=1}^{N} w_i g_{\varepsilon}''(\mathcal{J}(y_i,u,\xi_i)-t), \\
 \grad{\mathcal{L}}{u,u} & \approx \alpha \grad{\mathcal{P}(u)}{uu}, \\
 \grad{\mathcal{L}}{p_j,p_k}  & = 0,
\end{align*}
as well as symmetric terms.
For both notational and computational brevity, let us introduce
\begin{align*}
 H^{yt}_j & = -(1-\beta)^{-1} g_{\varepsilon}''(\mathcal{J}(y_j,u,\xi_j) - t) \grad{\mathcal{J}(y_j,u,\xi_j)}{y} \in Y \\
 H^{yy}_j & = (1-\beta)^{-1} \left[ g_{\varepsilon}''(\mathcal{J} - t) \grad{\mathcal{J}}{y} (\grad{\mathcal{J}}{y})^* + g_{\varepsilon}'(\mathcal{J} - t) \grad{\mathcal{J}}{yy} \right]_{j} \in \mathcal{L}(Y,Y) \\
 H^{tt} & = (1-\beta)^{-1} \sum_{i=1}^{N} w_i g_{\varepsilon}''(\mathcal{J}(y_i,u,\xi_i)-t) \in \mathbb{R}_+.
\end{align*}
The entire Newton KKT system can thus be written as follows,
\begin{align*}
& \begin{bmatrix}
 \mathrm{bdiag}(w_j H^{yy}_j) & 0 & \mathrm{bdiag}(w_j (\grad{c(y_j,u,\xi_j)}{y})^*) & w \odot H^{yt} \\
 0  & \alpha \grad{\mathcal{P}(u)}{uu} & [(\grad{c(y_k,u,\xi_k)}{u})^* w_k] & 0 \\
 \mathrm{bdiag}(w_j \grad{c(y_j,u,\xi_j)}{y}) & [w_j \grad{c(y_j,u,\xi_j)}{u}] & 0 & 0 \\
 (w \odot H^{yt})^* & 0 & 0 & H^{tt}
\end{bmatrix}
\begin{bmatrix}
 \delta_{\vec{y}} \\
 \delta_{u} \\
 \delta_{\vec{p}} \\
 \delta_t
\end{bmatrix}
= -\begin{bmatrix}
     \grad{\mathcal{L}}{\vec{y}} \\
     \grad{\mathcal{L}}{u} \\
     \grad{\mathcal{L}}{\vec{p}} \\
     \grad{\mathcal{L}}{t}
    \end{bmatrix},
\end{align*}
where $\mathrm{bdiag}(H_j)$ constructs a block-diagonal operator with $H_1,\ldots,H_N$ along the diagonal.

\subsection{Gauss-Newton system}
\label{s:GN}

Further simplifications of the Hessian can be introduced.
First, we can note that all terms corresponding to differentiating $\mathcal{L}$ in $y_j$ or $p_j$ first contain the quadrature weight $w_j$.
In higher dimensions $w_j$ may vary over many orders of magnitude, and keeping them in the Hessian may render it extremely ill-conditioned.
Instead, we divide the corresponding rows of the Hessian together with the right hand side (the derivatives of $\mathcal{L}$) by $w_j$ in our formulation directly.
This leads to a non-symmetric Gauss-Newton system, but a much lower condition number together with a potent preconditioner developed below allows one to use GMRES or BiCGStab with only a few iterations.

Second, in the case of a linear-quadratic control we have $\mathcal{P}(u) = \frac12 \langle u, M_u u\rangle$, and $c(y,u,\xi) = \hat c(y,\xi) + B u$. In a weakly nonlinear case we can consider an inexact Newton method, where $\mathcal{P}(u)$ and $c(y,u,\xi)$ are approximated in this form. This allows one to resolve \eqref{eq:dLdu} explicitly and plug the corresponding component $u = \mathrm{Proj}_{U_{ad}}((\alpha M_u)^{-1} (-B)^* \sum_{i=1}^{N} w_i p_i)$ back into the Lagrangian, reducing it to three variables $(\vec{y},\vec{p},t)$.
The reduced Gauss-Newton Hessian term appears:
\begin{equation}
\grad{\mathcal{L}}{p_j, p_k} = - w_j B P_r (\alpha M_u)^{-1} B^* w_k, \qquad P_r = \mathrm{Proj}'_{U_{ad}}\left((\alpha M_u)^{-1} (-B)^* \sum_{i=1}^{N} w_i p_i\right),
 \label{eq:dLdpdp}
\end{equation}
where $w_j$ in the leftmost position is pending to removal as described above,
and $\mathrm{Proj}'_{U_{ad}}$ is a semi-smooth derivative of the projector (e.g. for box constraints this is just an indicator of $U_{ad}$).

\begin{remark}
In addition to being quadratic, the control penalty is equivalent to a (weighted) squared $L_2$-norm in many cases.
This renders discretization of $M_u$ spectrally close to a diagonal matrix $\tilde M_u$; for example, one may use a standard lumping of the finite element mass matrix.
This makes the discretized Hessian \eqref{eq:dLdpdp} easy to assemble, e.g. sparse when $M_u$ is replaced by $\tilde M_u$.
The case of a $H^1$-norm control penalty is more limiting, and may require a matrix-free application of $M_u^{-1}$ using e.g. a multigrid method.
\end{remark}

The purpose of this elimination of $u$ becomes more apparent for the solution of the Gauss-Newton system in the TT format.
An Alternating Least Squares method (cf. Appendix~\ref{sec:ttalg}) tailored to Karush-Kuhn-Tucker (KKT) systems \cite{bdos-sb-2016} requires that all solution components are represented in the same TT decomposition \cite{dkos-eigb-2014}.
Since $\vec{y}$ and $\vec{p}$ have the same size, this holds naturally; the only additional variable $t$ is a single number that can be embedded into the same TT decomposition at a little cost.
In contrast, a (potentially large) component $u$ needs a nontrivial padding that may inflate the TT ranks and/or make the Hessian more ill-conditioned.

Finally, we obtain the following linear system on the solution increments:
\begin{equation}
\begin{bmatrix}
H^{yy} & A^* & H^{yt} \\
A & -\mathcal{B} \otimes W & 0 \\
(H^{yt}_w)^* & 0 & H^{tt}
\end{bmatrix}
\begin{bmatrix}
 \delta_{\vec{y}} \\
 \delta_{\vec{p}} \\
 \delta_t
\end{bmatrix} =
\begin{bmatrix}
F_{\vec{y}} \\
F_{\vec{p}} \\
F_t
\end{bmatrix},
\label{eq:GN}
\end{equation}
where $\mathcal{B} = B P_r (\alpha M_u)^{-1}B^*$, $H^{yy}$ and $A$ are (block) diagonal matrices with $H^{yy}_j$ and $A_j:=\grad{\hat c(y_j,\xi_j)}{y}$ on the diagonal,
$$
W = \begin{bmatrix} w_1 & \cdots & w_N \\ & \cdots & \\ w_1 & \cdots & w_N \end{bmatrix} \in \mathbb{R}^{N \times N}
$$
is a rank-1 matrix (after discretization), and
$$
H^{yt} = \begin{bmatrix}H^{yt}_1 \\ \vdots \\ H^{yt}_N\end{bmatrix}, \quad
H^{yt}_w = \begin{bmatrix} H^{yt}_1 w_1 \\ \vdots \\ H^{yt}_N w_N \end{bmatrix},
$$
whereas the right-hand side components are
\begin{align*}
 F_{\vec{y}} & = \frac{1}{1-\beta}\left[ g_{\varepsilon}' (\mathcal{J}(y_j,u,\xi_j) - t) \grad{\mathcal{J}(y_j,u,\xi_j)}{y} + \frac{1}{w_j M} \sum_{\ell=1}^{M} E_j'(\xi_\ell)  + (1-\beta)A_j^* p_j \right]_{j=1}^{N}, \\
 F_{\vec{p}} & = \left[\hat c(y_j,\xi_j) - \mathcal{B} \sum_{i=1}^N w_i p_i \right]_{j=1}^{N}, \\
 F_{t} & = 1 - \frac{1}{1-\beta}\left[ \sum_{i=1}^{N} w_i g_{\varepsilon}'(\mathcal{J}(y_i,u,\xi_i) - t) + \frac{1}{M}\sum_{\ell=1}^{M} e'(\xi_\ell)\right].
\end{align*}
Having solved \eqref{eq:GN},
we perform the Newton update similarly to Algorithm~\ref{alg:ttrisk}, by setting
$$
\vec{y}_{m+1} = \vec{y}_m + h\delta_{\vec{y}}, \quad \vec{p}_{m+1} = \vec{p}_m + h\delta_{\vec{p}}, \quad t_{m+1} = t_m + h \delta_t, \quad \mbox{and} \quad \varepsilon_{m+1} = \mu_{\varepsilon} \varepsilon_m,
$$
where the step size $h>0$ is chosen ensuring \eqref{eq:Arm1} and \eqref{eq:Arm2}.

\begin{remark}\label{rem:armijo}
 Note that \eqref{eq:Arm1}
 is not a proper line search condition and has no theoretical guarantees to lead to a convergent method. However, we have empirically observed in our numerical results that enforcing such a condition allowed our method to converge. A proper line search should be based for example on the Armijo condition applied either to the cost function (for the reduced space formulation in section~\ref{sec:reduced}) or a properly designed merit function (i.e. a weighted sum of the objective function and some norm of the constraint violation) for the full space formulation presented in this section.
\end{remark}

Observe that the system \eqref{eq:GN} can be ill-conditioned and thus requires a
good preconditioner to solve it efficiently. In what follows, we discuss a Schur 
complement-based preconditioner.

\subsection{Preconditioning}
\label{s:precond}

In this section we propose a matching strategy \cite{PSW} to approximate the Schur complement to the Gauss-Newton matrix~\eqref{eq:GN}.

First, since $\delta_t$ is a single number, we can compute the corresponding Schur complement matrix
\[
 \begin{bmatrix}
  H^{yy} - H^{yt}\frac{1}{H^{tt}} (H_w^{yt})^* & A^* \\
  A & -\mathcal{B} \otimes W
 \end{bmatrix},
\]
where we can denote $S^{yy}:=H^{yy} - H^{yt}\frac{1}{H^{tt}} (H_w^{yt})^*$ for brevity.
Since $A$ is a linearization of $c(y,u,\xi)$ (for example, the PDE operator), it is often invertible.
In this case, the Schur complement towards the $(1,2)$ block of this matrix reads
\[
 S = A^* + S^{yy} A^{-1} (\mathcal{B} \otimes W) = A^* A^{-1} A + S^{yy} A^{-1} (\mathcal{B} \otimes W).
\]
We propose a matching approximation consisting of three factors:
\begin{equation}
 \tilde S = (A^* + \eta S^{yy}) A^{-1} \left(A + \frac{1}{\eta}\mathcal{B} \otimes W \right),
\end{equation}
where $\eta = \sqrt{\|\mathcal{B} \otimes W\| / \|S^{yy}\|}$ is the scaling constant.
Note that
\[
 \|\tilde S - S\| = \|\eta S^{yy} + A^* A^{-1} \frac{1}{\eta}\mathcal{B} \otimes W\| = \mathcal{O}(\alpha^{-1/2}),
\]
which was shown \cite{bdos-sb-2016} to be small in norm compared to $S$ for both limits $\alpha \rightarrow \infty$ when $\|S\| = \mathcal{O}(1)$, and $\alpha\rightarrow 0$ when $\|S\| = \mathcal{O}(\alpha^{-1})$.

Ultimately, we obtain the following right preconditioner:
\begin{equation}\label{eq:prec}
 P = \begin{bmatrix}
      0 & \tilde S & H^{yt} \\
      A & -\mathcal{B}\otimes W & 0 \\
      0 & 0 & H^{tt}
     \end{bmatrix}.
\end{equation}
Note that this is a permuted block-triangular matrix, solving a linear system with which requires the solution of smaller systems with $H^{tt}, \tilde S$ and $A$.
In turn, the solution with $\tilde S$ requires the solution with $(A^* + \eta S^{yy})$ and $(A + \frac{1}{\eta}\mathcal{B} \otimes W)$.
If the constraints are defined by a PDE, $A,A^*$ and $g_{\varepsilon}'(\mathcal{J} - t) \grad{\mathcal{J}}{yy}$ (inside $H^{yy}$) are sparse matrices, whereas the remaining terms $g_{\varepsilon}''(\mathcal{J} - t) \grad{\mathcal{J}}{y} (\grad{\mathcal{J}}{y})^*$, $H^{yt}\frac{1}{H^{tt}} (H_w^{yt})^*$ and $W$ are low-rank matrices, and can be accounted for using the Sherman-Morrison formula.

\section{Numerical results}
\label{s:numerics}

In this section, we present various numerical examples to illustrate the efficiency of the proposed approach
in both the reduced and full-space formulations. This section in fact goes beyond the above theoretical
presentation in multiple ways. In section~\ref{s:1d}, we consider an optimal control problem in one spatial
dimension and random coefficient. We study the approximation error in CVAR$_\beta$ due to each of
the variables $\varepsilon$, $n_y$ (spatial discretization), $n_\xi,$, $d$, and 
TT truncation tolerance. We propose a
strategy to select these parameters by equidistribution of the total error. Control variate strategy is applied to
this problem in section~\ref{s:1d_cv}. Section~\ref{s:quantile} focuses on the impact of the quantile 
$\beta$ and the $\varepsilon$ reduction factor $\mu_\varepsilon$. In section~\ref{s:2d}, we consider the
two spatial dimension version of the
problem and carry out a comparison between the reduced and full space formulations. 
We conclude with a realistic problem in section~\ref{s:infection}, where we propose a 
risk-averse strategy for lockdown due to pandemics such as COVID-19.

The TTRISK (Algorithm~\ref{alg:ttrisk}) is implemented based on TT-Toolbox\footnote{\url{https://github.com/oseledets/TT-Toolbox}},
whereas for the TT-Cross (Algorithm~\ref{alg:tt-cross}) we use a rank-adaptive implementation \texttt{amen\_cross\_s} from TT-IRT\footnote{\url{https://github.com/dolgov/TT-IRT}}.
We run the computations using a default multithreading in Matlab R2019b that can spawn up to $10$ threads in BLAS on an Intel Xeon E5-2640 v4 CPU.

\subsection{Elliptic PDE with affine coefficient}\label{s:1d}
We first test the reduced formulation.
Consider a PDE in one space dimension with random coefficients $\kappa$.
\begin{align}\label{eq:PDE1}
-\frac{d}{dx} \kappa(x,\xi) \frac{dy}{dx} & = Bu, & x \in (0,1), \quad \xi \sim \mathcal{U}(-\sqrt{3},\sqrt{3})^d,\\
y(-1)=y(1)& =0, \nonumber
\end{align}
where $\mathcal{U}(-\sqrt{3},\sqrt{3})^d$ denotes the uniform distribution on $[-\sqrt{3},\sqrt{3}]^d$.
The Karhunen-Loeve (KL) expansion of $\kappa(x,\xi)$ truncated to $d$ terms,
\begin{equation}\label{eq:kle}
\kappa(x,\xi) = \kappa_0(x) + \sum_{k=1}^{d} \sqrt{\lambda_k} \kappa_k(x) \cdot \xi^{(k)},
\end{equation}
is defined by the mean $\kappa_0(x)=10$ and the eigenvalue decomposition
\begin{equation}\label{eq:kleeig}
\int C(x,x') \kappa_k(x') dx' = \lambda_k \kappa_k(x), \qquad \lambda_1 \ge \lambda_2 \ge \cdots \ge 0
\end{equation}
of the covariance operator with the function
$$
C(x,x') = \sigma^2 \exp\left(-\frac{(x-x')^2}{2\ell^2}\right), \qquad \ell=0.25,
$$
for $x , x' \in (0,1)$, and $\xi^{(k)} \in \mathcal{U}(-\sqrt{3},\sqrt{3})$, so that $\rho^{(k)}(\xi^{(k)}) = 1/(2\sqrt{3})$.

We use a misfit objective function,
\begin{equation}\label{eq:J1d}
 \mathcal{J}(y,u;\xi) = \frac{1}{2}\|y(u; x,\xi) - y_d(x)\|_{L_2(0,1)}^2,
\end{equation}
with the desired state $y_d(x) \equiv 1$.
The control $u$ is defined on a subdomain $(0.25, 0.75)$, and $B$ is the identity insertion operator:
$$
Bu(x) = \left\{\begin{array}{ll}u(x), & x \in (0.25, 0.75), \\ 0, & \mbox{otherwise.}\end{array}\right.
$$
The PDE is discretised using piecewise linear finite elements on a uniform grid with $n_y$ points $0=x_1<\cdots<x_{n_y}=1$, with the coefficient $\kappa(x,\xi)$ and the control $u(x)$ discretised by collocation at the midpoints $\{x_{i+1/2}\}$.
The random variables $\xi^{(1)},\ldots,\xi^{(d)}$ are discretised by collocation at $n_{\xi}$ Gauss-Legendre points on the interval $(-\sqrt{3},\sqrt{3})$.

We aim at estimating the total cost-error scaling.
However, since the computation depends on a number of approximation parameters ($\varepsilon$, $n_y$, $n_{\xi}$, $d$ and the TT truncation tolerance), those need to be selected judiciously to obtain the optimal total complexity.
Next, we estimate the error (in CVaR) contributed by each of the parameters.
This will allow us to select the parameter values by equalizing their corresponding error estimates.

The (relative) CVaR error at given parameters is defined as
\begin{equation}\label{eq:cvar-err}
\mbox{err(CVaR)} = \frac{\mathcal{R}_{t_{m},\beta,N}^{\varepsilon} - \mathcal{R}_*}{\mathcal{R}_*},
\end{equation}
where
$t_{m}$ is the output of Algorithm~\ref{alg:ttrisk}, and $\mathcal{R}_*$ is
the reference solution computed at the finest parameters $\varepsilon=3\cdot 10^{-4}$, tolerance$=10^{-5}$, $n_y=1025$, $n_{\xi}=33$, $d=10$.
We take the control regularization parameter $\alpha=10^{-6}$ and the confidence threshold $\beta=0.5$.
In the following figures, we vary those parameters one by one, keeping the other fixed to those reference values.

%

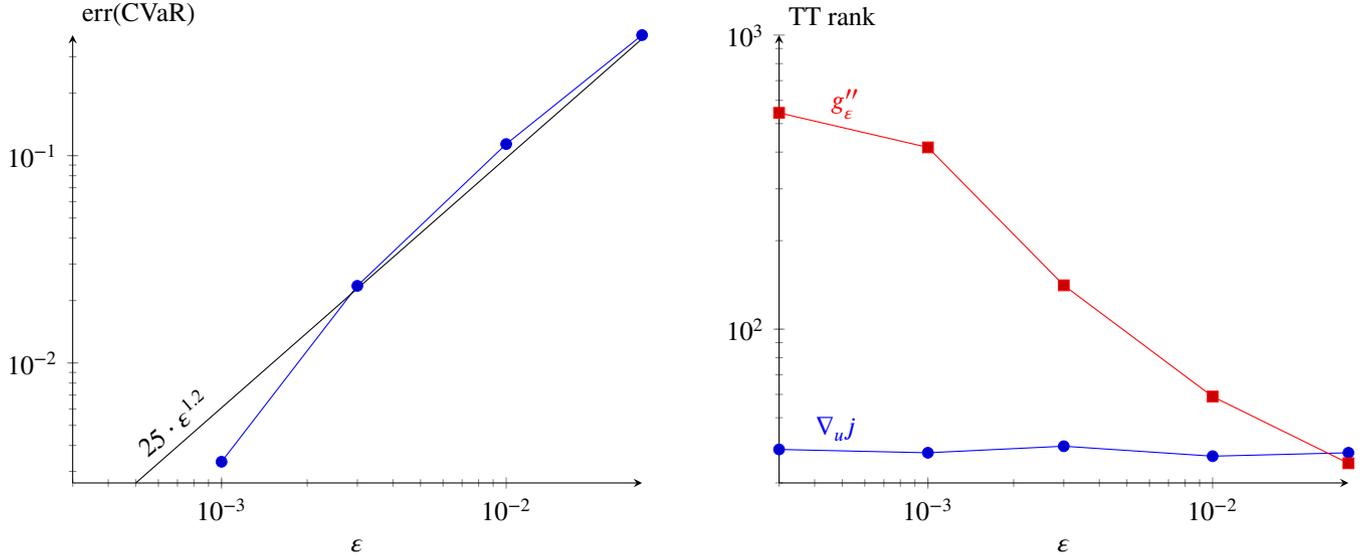
\begin{figure}[htb]
\beginpgfgraphicnamed{CVaR_ML_Rev4-fig-err-eps1}
\begin{tikzpicture}
\begin{axis}[%
 width=0.42\linewidth,
 height=0.33\linewidth,
 xlabel=$\varepsilon$,
 ylabel=err(CVaR),
 xmode=log,
 ymode=log,
 xmin=3e-4,xmax=3e-2,
 ]
 \addplot+[] coordinates{
 (3e-2, 3.8133e-01)
 (1e-2, 1.1367e-01)
 (3e-3, 2.3544e-02)
 (1e-3, 3.3372e-03)
 };
 \addplot+[domain=5e-4:3e-2,no marks,black] {x^(1/0.83)*25} node[pos=0.1,anchor=west,above,rotate=40] {$25\cdot \varepsilon^{1.2}$};
\end{axis}
\end{tikzpicture}
\endpgfgraphicnamed
\hfill
\beginpgfgraphicnamed{CVaR_ML_Rev4-fig-err-eps2}
\begin{tikzpicture}
\begin{axis}[%
 width=0.42\linewidth,
 height=0.33\linewidth,
 xlabel=$\varepsilon$,
 ylabel=TT rank,
 xmode=log,
 ymode=log,
 xmin=3e-4,xmax=3e-2,
 ymin=30,ymax=1e3,
 ]
 \addplot+[] coordinates{
 (3e-2, 38)
 (1e-2, 37)
 (3e-3, 40)
 (1e-3, 38)
 (3e-4, 39)
 } node[pos=0.9,anchor=west,above] {$\grad{j}{u}$};
 \addplot+[] coordinates{
 (3e-2, 35 )
 (1e-2, 59 )
 (3e-3, 141)
 (1e-3, 415)
 (3e-4, 544)
 } node[pos=0.9,anchor=west,above] {$g''_{\varepsilon}$};
\end{axis}
\end{tikzpicture}
\endpgfgraphicnamed
\caption{CVaR error \eqref{eq:cvar-err} and TT ranks depending on the CVaR smoothing parameter. Other parameters: tolerance$=10^{-5}$, $n_y=1025$, $n_{\xi}=33$, $d=10$, $\alpha=10^{-6}$, $\beta=0.5$. \label{fig:err-eps}}
\end{figure}

Figure~\ref{fig:err-eps} shows that the error in CVaR depends almost linearly on $\varepsilon$ (in fact the decay is slightly faster, which may be due to particular symmetries in the solution).
This is consistent with \cite[Lemma 3.4.2]{MaterThesisMae} where a linear theoretical convergence was established.
However, the TT ranks of the logistic function derivatives grow also linearly with $1/\varepsilon$.
This will eventually lead to noticeable computing costs as $\varepsilon$ decreases.
It is thus crucial that the cross approximation of $g''_{\varepsilon}$ uses a precomputed TT approximation of $j(u,\xi)$ instead of original PDE solutions.
In turn, the TT ranks of $j(u,\xi)$ and $\grad{j(u,\xi)}{u}$ stay almost constant near $40$, and hence the number of PDE solutions is almost independent
of $\varepsilon$ in the considered range.


\begin{figure}[htb]
\beginpgfgraphicnamed{CVaR_ML_Rev4-fig-err-nxi1}
\begin{tikzpicture}
\begin{axis}[%
 width=0.41\linewidth,
 height=0.33\linewidth,
 xlabel=$n_{\xi}$,
 ylabel=err(CVaR),
 xmode=normal,
 ymode=log,
 ]
 \addplot+[] coordinates{
 (3, 2.8998e-03)
 (4, 8.3174e-04)
 (5, 3.0911e-04)
 (7, 7.5365e-05)
 };
 \addplot+[domain=3:7,no marks,black] {exp(-x)*1e-1} node [pos=0.4,anchor=west,above,rotate=-38] {$10^{-1}\cdot \exp(-n_{\xi})$};
\end{axis}
\end{tikzpicture}
\endpgfgraphicnamed
\hfill
\beginpgfgraphicnamed{CVaR_ML_Rev4-fig-err-nxi2}
\begin{tikzpicture}
\begin{axis}[%
 width=0.41\linewidth,
 height=0.33\linewidth,
 xlabel=$n_{\xi}$,
 ylabel=TT rank,
 xmode=normal,
 ymode=log,
 ymin=30,ymax=5e2,
 ]
 \addplot+[] coordinates{
 (3, 38)
 (4, 37)
 (5, 37)
 (7, 39)
 } node[pos=0.9,anchor=east,above] {$\grad{j}{u}$};
 \addplot+[] coordinates{
 (3, 71 )
 (4, 154)
 (5, 289)
 (7, 486)
 } node[pos=0.9,anchor=east,below] {$g''_{\varepsilon}$};
\end{axis}
\end{tikzpicture}
\endpgfgraphicnamed
\caption{CVaR error \eqref{eq:cvar-err} and TT ranks depending on the random variable discretization. Other parameters: $\varepsilon=3\cdot 10^{-4}$, tolerance$=10^{-5}$, $n_y=1025$, $d=10$, $\alpha=10^{-6}$, $\beta=0.5$. \label{fig:err-nxi}}
\end{figure}
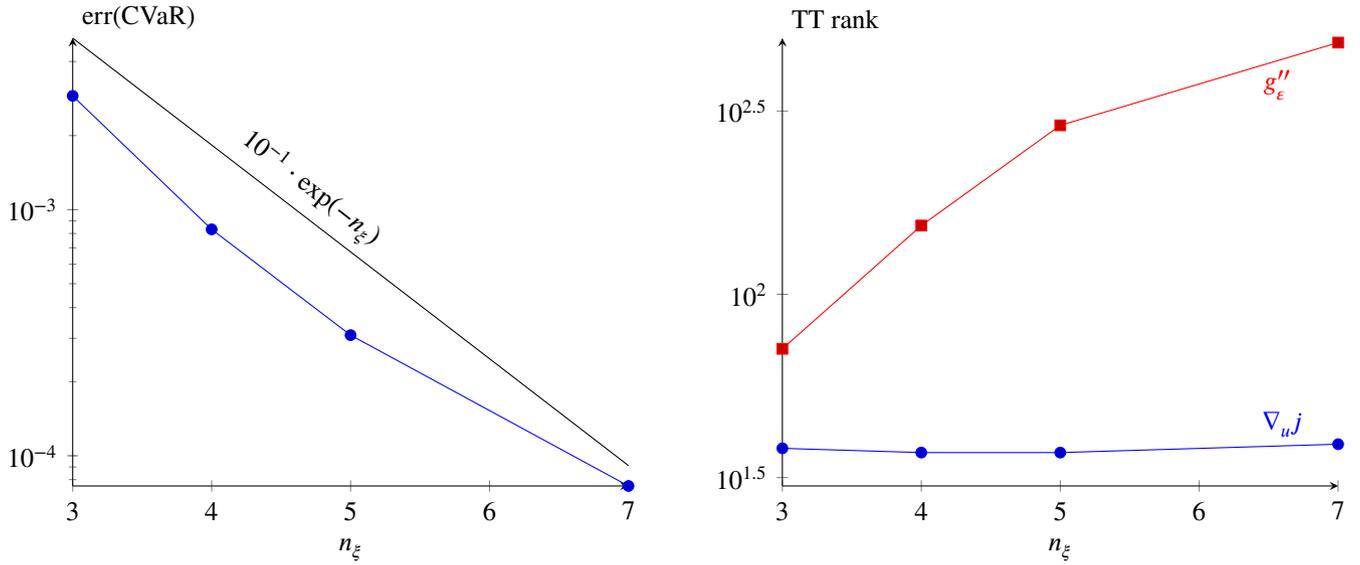

In Figure~\ref{fig:err-nxi} we vary the number of quadrature points introduced in each of the random variables $\xi$.
As expected from many previous works (see e.g. \cite{KS11,schwab-analyticity-2011,Herrmann2020}),
the approximation converges exponentially in $n_{\xi}$.
The TT ranks stabilize towards the value prescribed by $\varepsilon$ in Figure~\ref{fig:err-eps}.





\begin{figure}[htb]
\beginpgfgraphicnamed{CVaR_ML_Rev4-fig-err-tol1}
\begin{tikzpicture}
\begin{axis}[%
 width=0.35\linewidth,
 height=0.30\linewidth,
 xlabel=tolerance,
 ylabel=err(CVaR),
 xmode=log,
 ymode=log,
 ]
 \addplot+[] coordinates{
 (1e-2, 4.1590e+00)
 (3e-3, 1.8913e-02)
 (1e-3, 3.6471e-03)
 (3e-4, 1.0007e-03)
 (1e-4, 3.1684e-04)
 };
 \addplot+[domain=1e-4:1e-2,no marks,black] {x*4} node [pos=0.9,anchor=east,below,rotate=22] {$4\cdot\mbox{tol}$};
\end{axis}
\end{tikzpicture}
\endpgfgraphicnamed
\hfill
\beginpgfgraphicnamed{CVaR_ML_Rev4-fig-err-tol2}
\begin{tikzpicture}
\begin{axis}[%
 width=0.35\linewidth,
 height=0.30\linewidth,
 xlabel=tolerance,
 ylabel=TT rank,
 xmode=log,
 ymode=log,
 ]
 \addplot+[] coordinates{
 (1e-2, 9 )
 (3e-3, 12)
 (1e-3, 15)
 (3e-4, 20)
 (1e-4, 24)
 } node[pos=0.9,anchor=east,above] {$\grad{j}{u}$};
 \addplot+[] coordinates{
 (1e-2, 57 )
 (3e-3, 589)
 (1e-3, 752)
 (3e-4, 396)
 (1e-4, 422)
 } node[pos=0.9,anchor=east,below] {$g''_{\varepsilon}$};
\end{axis}
\end{tikzpicture}
\endpgfgraphicnamed
\caption{CVaR error \eqref{eq:cvar-err} and TT ranks depending on the TT truncation tolerance. Other parameters: $\varepsilon=3\cdot 10^{-4}$, $n_y=1025$, $n_{\xi}=33$, $d=10$, $\alpha=10^{-6}$, $\beta=0.5$. \label{fig:err-tol}}
\end{figure}
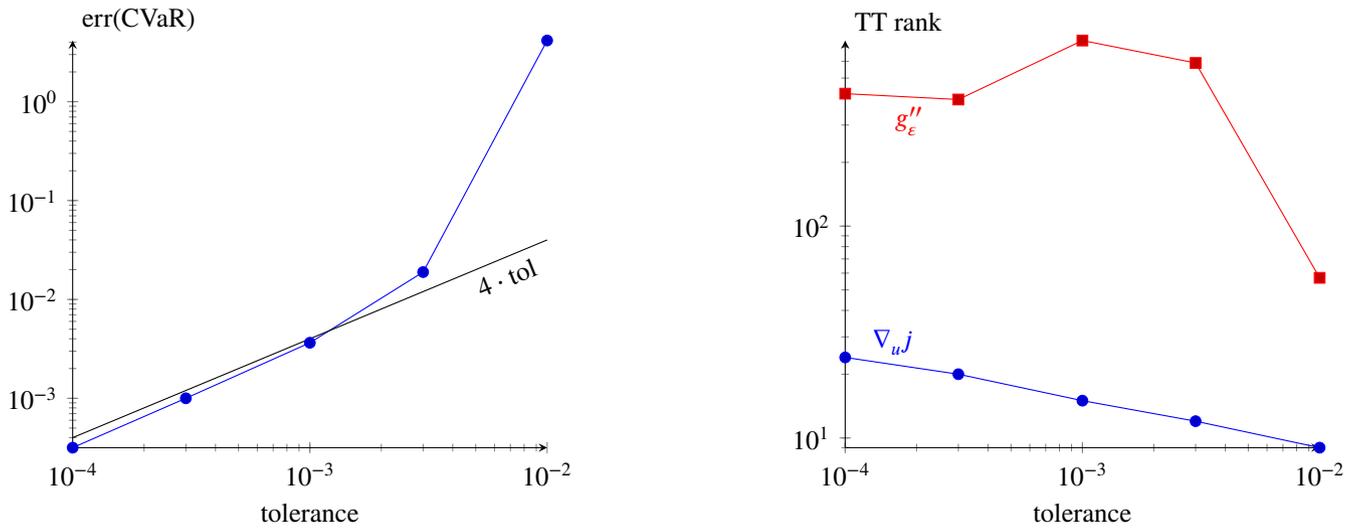

Figure~\ref{fig:err-tol} benchmarks the scheme against the relative error tolerance used for stopping both the TT-Cross Algorithm~\ref{alg:tt-cross} and TTRISK Algorithm~\ref{alg:ttrisk}, as well as for the TT approximation.
The convergence is linear except for very large values of the threshold, when the TT-Cross may stop prematurely after an accidental drop of the iteration increment below the threshold.
The TT ranks grow logarithmically or even slower with $1/\mbox{tol}$, which is the enabling observation for many applications of tensor methods.


\begin{figure}[htb]
\beginpgfgraphicnamed{CVaR_ML_Rev4-fig-err-ny1}
\begin{tikzpicture}
\begin{axis}[%
 width=0.40\linewidth,
 height=0.33\linewidth,
 xlabel=$n_{y}-1$,
 ylabel=err(CVaR),
 xmode=log,
 ymode=log,
 xtick={64,128,256,512},
 xticklabels={64,128,256,512},
 ]
 \addplot+[] coordinates{
 (65 -1, 1.6355e-03)
 (129-1, 3.4115e-04)
 (257-1, 2.4586e-04)
 (513-1, 2.7317e-05)
 };
 \addplot+[domain=64:512,no marks,black] {x^(-2)*10} node [pos=0.3,anchor=west,above,rotate=-38] {$10 (n_y-1)^{-2}$};
\end{axis}
\end{tikzpicture}
\endpgfgraphicnamed
\hfill
\beginpgfgraphicnamed{CVaR_ML_Rev4-fig-err-ny2}
\begin{tikzpicture}
\begin{axis}[%
 width=0.40\linewidth,
 height=0.33\linewidth,
 xlabel=$n_{y}-1$,
 ylabel=TT rank,
 xmode=log,
 ymode=log,
 ymin=30,ymax=1e3,
 xtick={64,128,256,512},
 xticklabels={64,128,256,512},
 ]
 \addplot+[] coordinates{
 (65 -1, 36)
 (129-1, 36)
 (257-1, 37)
 (513-1, 37)
 } node[pos=0.9,anchor=east,above] {$\grad{j}{u}$};
 \addplot+[] coordinates{
 (65 -1, 464)
 (129-1, 463)
 (257-1, 869)
 (513-1, 798)
 } node[pos=0.9,anchor=east,below] {$g''_{\varepsilon}$};
\end{axis}
\end{tikzpicture}
\endpgfgraphicnamed
\caption{CVaR error \eqref{eq:cvar-err} and TT ranks depending on the spatial discretization. Other parameters: $\varepsilon=3\cdot 10^{-4}$, tolerance$=10^{-5}$, $n_{\xi}=33$, $d=10$, $\alpha=10^{-6}$, $\beta=0.5$. \label{fig:err-ny}}
\end{figure}
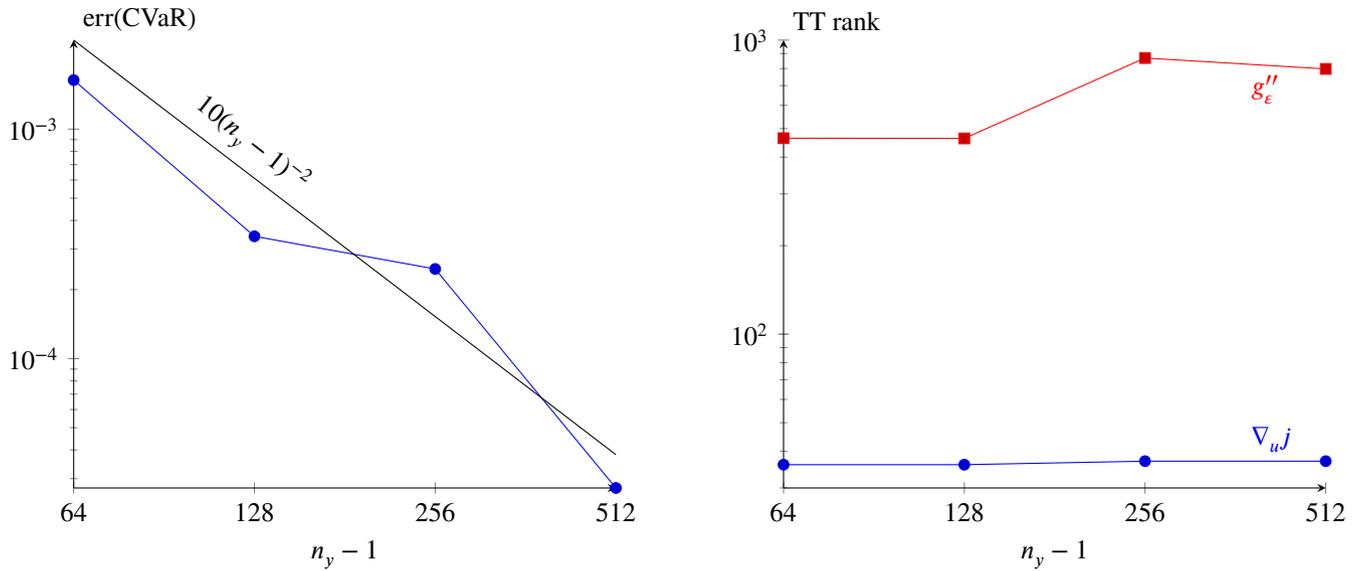

In Figure~\ref{fig:err-ny} we vary the number of finite elements in the spatial discretization of \eqref{eq:PDE1}.
As expected from the second order of consistency of the continuous linear elements, the CVaR error converges quadratically with $n_y$.
The TT ranks stay almost constant, which shows that even the coarsest grid resembles enough qualitative features of the solution.


\begin{figure}[htb]
\beginpgfgraphicnamed{CVaR_ML_Rev4-fig-err-d1}
\begin{tikzpicture}
\begin{axis}[%
 width=0.40\linewidth,
 height=0.33\linewidth,
 xlabel=$d$,
 ylabel=err(CVaR),
 xmode=normal,
 ymode=log,
 ]
 \addplot+[] coordinates{
 (2, 1.5314e-02)
 (3, 9.4047e-04)
 (4, 9.2442e-04)
 (5, 4.7411e-05)
 (6, 4.0004e-05)
 };
 \addplot+[domain=2:6,no marks,black] {exp(-2*x)} node [pos=0.8,anchor=east,below,rotate=-38] {$\exp(-2d)$};
\end{axis}
\end{tikzpicture}
\endpgfgraphicnamed
\hfill
\beginpgfgraphicnamed{CVaR_ML_Rev4-fig-err-d2}
\begin{tikzpicture}
\begin{axis}[%
 width=0.40\linewidth,
 height=0.33\linewidth,
 xlabel=$d$,
 ylabel=TT rank,
 xmode=normal,
 ymode=log,
 ]
 \addplot+[] coordinates{
 (2, 6 )
 (3, 12)
 (4, 19)
 (5, 24)
 (6, 31)
 } node[pos=0.9,anchor=east,above] {$\grad{j}{u}$};
 \addplot+[] coordinates{
 (2, 11 )
 (3, 36 )
 (4, 165)
 (5, 190)
 (6, 199)
 } node[pos=0.9,anchor=east,below] {$g''_{\varepsilon}$};
\end{axis}
\end{tikzpicture}
\endpgfgraphicnamed
\caption{CVaR error \eqref{eq:cvar-err} and TT ranks depending on the number of terms in the KL expansion. Other parameters: $\varepsilon=3\cdot 10^{-4}$, tolerance$=10^{-5}$, $n_y=1025$, $n_{\xi}=33$, $\alpha=10^{-6}$, $\beta=0.5$. \label{fig:err-d}}
\end{figure}
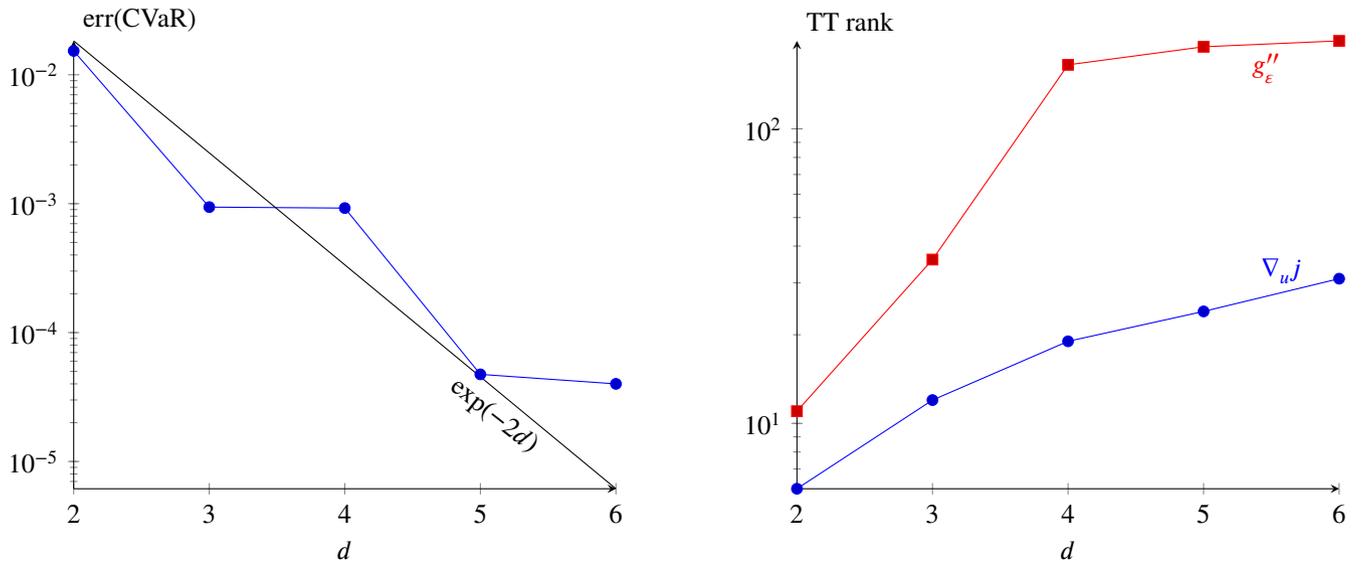

Lastly in this series, Figure~\ref{fig:err-d} varies the dimension of the random vector $\xi$, that is, the number of terms in the KL expansion.
For the Gaussian covariance matrix of $\kappa$ we observe the expected exponential convergence in $d$,
and the stabilization of the TT ranks as long as the contribution of the latter random variables ($\xi^{(k)},\ldots,\xi^{(d)}$ for $k\ge k_0$ with some $k_0>1$) becomes negligible compared to the (fixed)
TT truncation threshold.

Equipped with the individual error estimates,
we can return to estimating the total error-cost scaling.
We choose all approximation parameters ($n_y,n_{\xi},\varepsilon,d$ and tolerance) such that errors predicted using the rules fitted in Figures~\ref{fig:err-eps}--\ref{fig:err-d} are equal.
Expanding the approximate solution as a Taylor series around the exact solution,
we obtain that up to second-order terms, the total error is less than $5$ times the individual error.
However, instead of varying directly the total error and calculating all parameters accordingly, it is more convenient to vary the most discrete parameter, to estimate the corresponding error contribution, and to calculate the remaining $4$ parameters using the inverse error prediction rules.
The most discrete parameter is $n_y$
(the spatial grid size), since the number of grid intervals is restricted to a multiple of 4 to ensure that the control subdomain is aligned to all grids considered.
For each $n_y$, this gives the estimated part of the error
$$
\frac{\mbox{error}_{\mbox{total}}}{5} = 10(n_y-1)^{-2},
$$
and the other parameters are selected as
\begin{align*}
\varepsilon &= \left(\frac{\mbox{error}_{\mbox{total}}}{25\cdot 5}\right)^{0.83}, & n_{\xi} &= \left\lceil -\log\left(10 \frac{\mbox{error}_{\mbox{total}}}{5}\right) \right\rceil, \\
\mbox{tol} &= \frac{\mbox{error}_{\mbox{total}}}{4 \cdot 5}, & d &= \left\lceil-\frac{1}{2}\log\left(\frac{\mbox{error}_{\mbox{total}}}{5}\right) \right\rceil.
\end{align*}


\begin{figure}[htb]
\beginpgfgraphicnamed{CVaR_ML_Rev4-fig-err-total1}
\begin{tikzpicture}
\begin{axis}[%
 width=0.42\linewidth,
 height=0.33\linewidth,
 xlabel=err(CVaR),
 xmode=log,
 ymode=log,
 ]
 \addplot+[] coordinates{
 (1.4540e-02, 12260 /1000)
 (1.6896e-03, 59995 /1000)
 (7.3791e-04, 68015 /1000)
 (3.8528e-04, 161395/1000)
 (1.6171e-04, 243135/1000)
 } node[pos=0.3,anchor=east,below,rotate=-37] {thousands PDE solves};
 \addplot+[] coordinates{
 (1.4540e-02, 17)
 (1.6896e-03, 75)
 (7.3791e-04, 99)
 (3.8528e-04, 276)
 (1.6171e-04, 293)
 } node[pos=0.2,anchor=east,above,rotate=-36] {TT rank($g''_{\varepsilon}$)};
 \addplot+[mark=diamond*] coordinates{
 (1.4540e-02, 7 )
 (1.6896e-03, 11)
 (7.3791e-04, 12)
 (3.8528e-04, 17)
 (1.6171e-04, 18)
 } node[pos=0.8,anchor=east,above,rotate=-10] {TT rank($j$)};
 \addplot+[domain=1e-4:1e-2,no marks,black] {1.2*x^(-0.6)} node [pos=0.2,anchor=east,below,rotate=-35] {$\mbox{err}^{-0.6}$};
\end{axis}
\end{tikzpicture}
\endpgfgraphicnamed
\hfill
\beginpgfgraphicnamed{CVaR_ML_Rev4-fig-err-total2}
\begin{tikzpicture}
\begin{axis}[%
 width=0.42\linewidth,
 height=0.33\linewidth,
 xlabel=err(CVaR),
 xmode=log,
 ymode=log,
 ]
 \addplot+[] coordinates{
 (1.4540e-02,     0.7552)
 (1.6896e-03,     4.2467)
 (7.3791e-04,     6.4615)
 (3.8528e-04,    81.7455)
 (1.6171e-04,   179.8575)
 } node[pos=0.3,below,rotate=-27] {CPU time (sec.)};
 \addplot+[domain=1e-4:1e-2,no marks,black] {0.015*x^(-1)} node [pos=0.2,anchor=east,below,rotate=-35] {$\mbox{err}^{-1}$};
\end{axis}
\end{tikzpicture}
\endpgfgraphicnamed
\caption{Total CVaR error \eqref{eq:cvar-err} and computing cost with parameters $d,n_y,n_{\xi},\mbox{tol},\varepsilon$ chosen to equilibrate their individual error contributions. \label{fig:err-total}}
\end{figure}
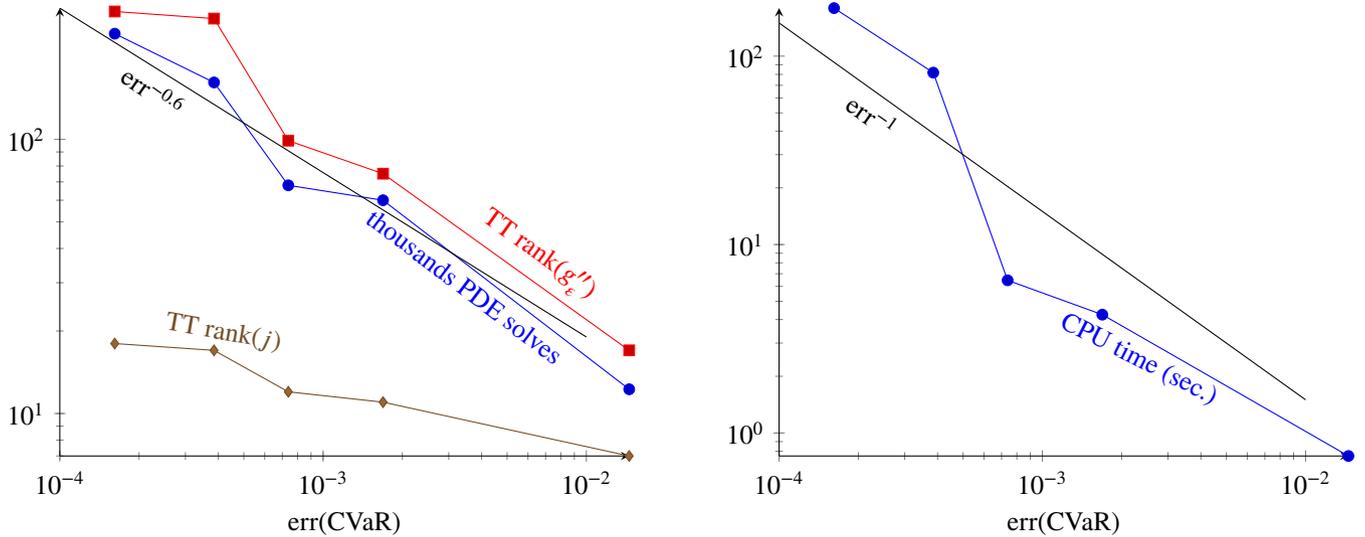

In Figure~\ref{fig:err-total}
we show the TT ranks, the corresponding numbers of PDE solutions required for the TT-Cross approximation of $j(u,\xi)$ and $\grad{j(u,\xi)}{u}$,
as well as the total computing time,
as functions of the \emph{actual} total error for $n_y=33,65,129,193$ and $257$, and the other parameters chosen as above.
We observe that both complexity indicators depend linearly on $1/\mbox{error}$ or slower.
This is to be compared with the solution of a deterministic PDE \eqref{eq:PDE1} with linear finite elements,
which provide an $\mbox{error} = O(n_y^{-2})$ with the computing $\mbox{cost} = O(n_y)$,
that is, the deterministic problem scales as $\mbox{cost} = \mbox{error}^{-1/2}$.
We see that the TT solution of a high-dimensional stochastic problem contributes the same amount of complexity.

\subsection{Control variate correction}\label{s:1d_cv}

From Figure~\ref{fig:err-eps} we notice that $\varepsilon$ is the slowest in terms of convergence.
In this experiment we add the Monte Carlo control variate correction \eqref{eq:MCgradu}, \eqref{eq:MCgradt}
to the gradient of the CVaR cost function during the optimization, and to the CVaR computation \eqref{eq:MCred} in the end.
Note that the first terms in \eqref{eq:MCred}--\eqref{eq:MCgradt} are deterministic.
Therefore, the standard deviations of \eqref{eq:MCred}--\eqref{eq:MCgradt} are equal to the standard deviations of the correction terms only, and can be seen as errors in the Euclidean norm of the underlying probability space.
To estimate these errors numerically, we run the algorithm (using some $\tilde M$ samples of $\xi$ to compute the averages in \eqref{eq:MCred}--\eqref{eq:MCgradt}) $16$ times, which gives us $16$ iid samples of \eqref{eq:MCred}--\eqref{eq:MCgradt},
and compute empirical standard deviations of those.
The ultimate correction is computed as another average of those $16$ corrections,
so we let $M=16\tilde M$ denote the total number of samples of $\xi$ from all runs of the algorithm.
The standard deviation of this ultimate result can be estimated as $1/4$-th of the empirical standard deviation computed above.

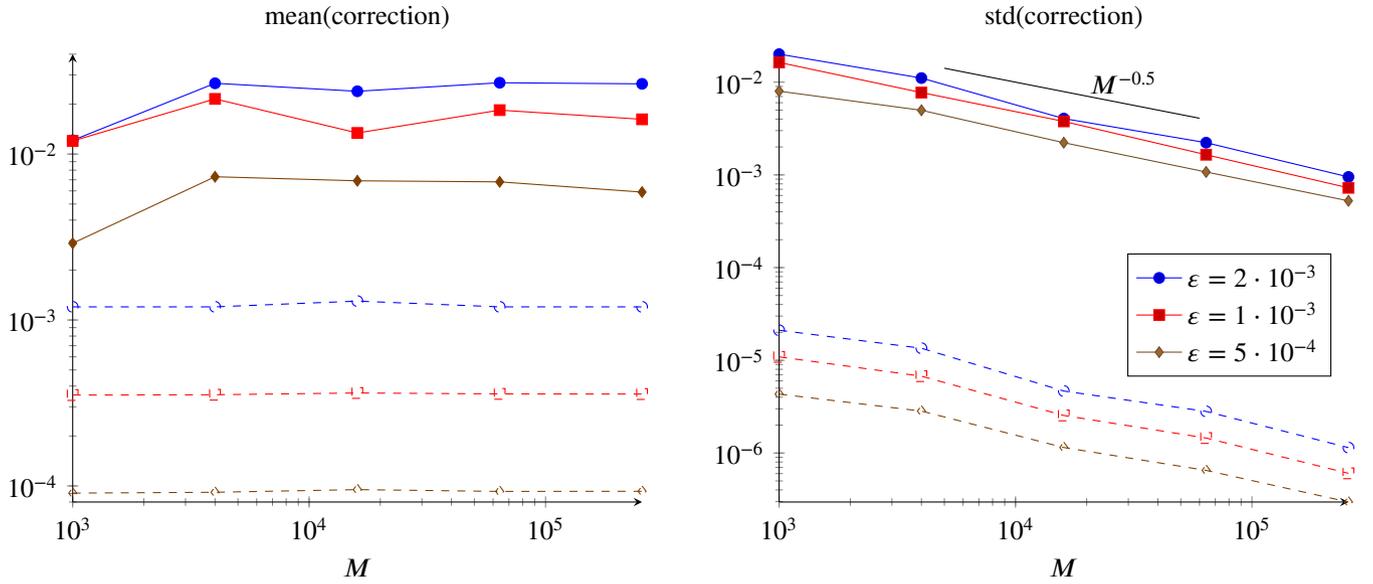
\begin{figure}[htb]
\beginpgfgraphicnamed{CVaR_ML_Rev4-fig-MLMC1}
\begin{tikzpicture}
\begin{axis}[%
 width=0.42\linewidth,
 height=0.33\linewidth,
 xlabel=$M$,
 title=mean(correction),
 xmode=log,
 ymode=log,
 legend style={at={(0.97,0.03)},anchor=south east},
 ymin=8e-5,ymax=4e-2,
 ]
 \addplot+[blue,solid,mark=*,mark options={style=blue}] coordinates{
 ( 1e3, 0.0121)
 ( 4e3, 0.0267 )
 (16e3, 0.0239 )
 (64e3, 0.0269 )
(256e3, 0.0265 )
 };
 \addplot+[red,solid,mark=square*,mark options={style=red}] coordinates{
 ( 1e3, 0.0120)
 ( 4e3, 0.0215 )
 (16e3, 0.0134 )
 (64e3, 0.0184 )
(256e3, 0.0162 )
 };
 \addplot+[orange!50!black,solid,mark=diamond*,mark options={style=orange!50!black}] coordinates{
 ( 1e3, 0.0029)
 ( 4e3, 0.0073 )
 (16e3, 0.0069 )
 (64e3, 0.0068 )
(256e3, 0.0059 )
 };
 \addplot+[blue,dashed,mark=o,mark options={style=blue}] coordinates{
 ( 1e3, 0.0012)
 ( 4e3, 0.0012)
 (16e3, 0.0013)
 (64e3, 0.0012)
(256e3, 0.0012)
 };
 \addplot+[red,dashed,mark=square,mark options={style=red}] coordinates{
 ( 1e3, 3.5303e-04)
 ( 4e3, 3.5451e-04)
 (16e3, 3.6330e-04)
 (64e3, 3.5867e-04)
(256e3, 3.5791e-04)
 };
 \addplot+[orange!50!black,dashed,mark=diamond,mark options={style=orange!50!black}] coordinates{
 ( 1e3, 9.0530e-05)
 ( 4e3, 9.1581e-05)
 (16e3, 9.5029e-05)
 (64e3, 9.2677e-05)
(256e3, 9.2846e-05)
 };
\end{axis}
\end{tikzpicture}
\endpgfgraphicnamed
\hfill
\beginpgfgraphicnamed{CVaR_ML_Rev4-fig-MLMC2}
\begin{tikzpicture}
\begin{axis}[%
 width=0.42\linewidth,
 height=0.33\linewidth,
 xlabel=$M$,
 title=std(correction),
 xmode=log,
 ymode=log,
 legend style={at={(0.97,0.28)},anchor=south east},
 ]
 \addplot+[] coordinates{
 ( 1e3, 0.0802/4)
 ( 4e3, 0.0443/4)
 (16e3, 0.0162/4)
 (64e3, 0.0089/4)
(256e3, 0.0038/4 )
 }; \addlegendentry{$\varepsilon=2\cdot 10^{-3}$};
 \addplot+[] coordinates{
 ( 1e3, 0.0653/4)
 ( 4e3, 0.0309/4)
 (16e3, 0.0151/4)
 (64e3, 0.0066/4)
(256e3, 0.0029/4 )
 }; \addlegendentry{$\varepsilon=1\cdot 10^{-3}$};
 \addplot+[mark=diamond*] coordinates{
 ( 1e3, 0.0320/4)
 ( 4e3, 0.0199/4)
 (16e3, 0.0089/4)
 (64e3, 0.0043/4)
(256e3, 0.0021/4 )
 }; \addlegendentry{$\varepsilon=5\cdot 10^{-4}$};
 \addplot+[blue,dashed,mark=o,mark options={style=blue}] coordinates{
 ( 1e3, 8.4412e-05/4)
 ( 4e3, 5.3973e-05/4)
 (16e3, 1.8651e-05/4)
 (64e3, 1.1243e-05/4)
 (256e3,4.5752e-06/4)
 };
 \addplot+[red,dashed,mark=square,mark options={style=red}] coordinates{
 ( 1e3, 4.3541e-05/4)
 ( 4e3, 2.6976e-05/4)
 (16e3, 1.0172e-05/4)
 (64e3, 5.8075e-06/4)
 (256e3,2.4151e-06/4)
 };
 \addplot+[orange!50!black,dashed,mark=diamond,mark options={style=orange!50!black}] coordinates{
 ( 1e3, 1.7292e-05/4)
 ( 4e3, 1.1335e-05/4)
 (16e3, 4.5888e-06/4)
 (64e3, 2.6000e-06/4)
 (256e3,1.1848e-06/4)
 };

 \addplot[black,solid,no marks,domain=5e3:6e4] {x^(-0.5)} node[pos=0.7,above] {$M^{-0.5}$};
\end{axis}
\end{tikzpicture}
\endpgfgraphicnamed
\caption{Solid lines: average (left) and estimated standard deviation (right) of the correction to the cost gradient \eqref{eq:MCgradu},\eqref{eq:MCgradt} depending on the number of Monte Carlo samples $M$. Dashed lines: average and standard deviation of the correction to $\mathcal{R}_{t}$ \eqref{eq:MCred}. Spatial discretization $n_y=129$, other parameters are chosen to equilibrate their individual error contributions. \label{fig:MLMC}}
\end{figure}

In Figure~\ref{fig:MLMC} we show both the ultimate average corrections and standard deviations of the corrections estimated as above.
We see that the standard deviations decay with the law of large numbers, as expected.
Moreover, both means and standard deviations decrease linearly with $\varepsilon$.

In multilevel Monte Carlo,
the estimated standard deviation can also be used to adapt the number of Monte Carlo samples towards the desired error threshold by extrapolating the law of large  numbers \cite{Scheichl-mlqmc-lognorm-2017}.
However, the convergence (and error estimation) of the TT approximation can be more complicated.
Therefore, we suggest to decrease $\varepsilon$ until the TT ranks are still manageable, and then compute the control variate correction to both estimate and improve the error.
This also allows us to decouple the TT and Monte Carlo steps, or to reuse previously computed TT approximations of e.g. the forward model solution.

\subsection{Effect of the quantile}\label{s:quantile}

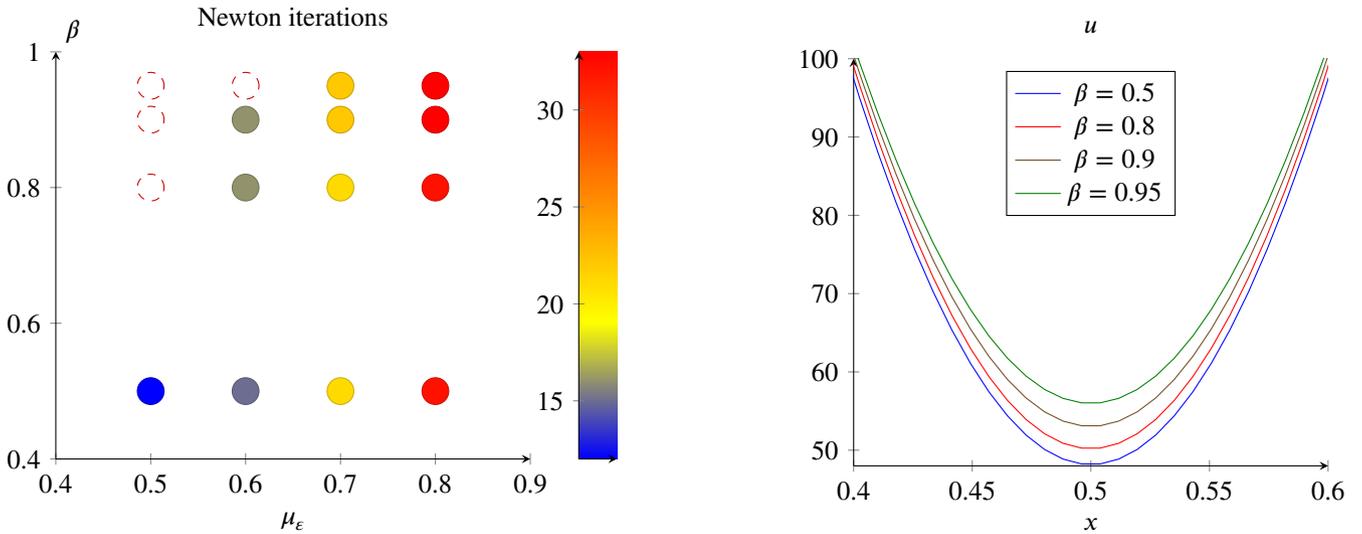
\begin{figure}[htb]
\beginpgfgraphicnamed{CVaR_ML_Rev4-fig-beta1}
\begin{tikzpicture}
 \begin{axis}[%
 width=0.35\linewidth,
 height=0.30\linewidth,
 view={0}{90},
 ylabel=$\beta$,
 xlabel=$\mu_{\varepsilon}$,
 colorbar,
 xmin=0.4,xmax=0.9,
 ymin=0.4,ymax=1.0,
 title={Newton iterations},
 ]
\addplot3[only marks,scatter,mark size=5] coordinates{
(0.5, 0.5 ,          12 )
(0.6, 0.5 ,          15 )
(0.7, 0.5 ,          21 )
(0.8, 0.5 ,          32 )
(0.5, 0.8 ,          nan)
(0.6, 0.8 ,          16 )
(0.7, 0.8 ,          21 )
(0.8, 0.8 ,          32 )
(0.5, 0.9 ,          nan)
(0.6, 0.9 ,          16 )
(0.7, 0.9 ,          22 )
(0.8, 0.9 ,          33 )
(0.5, 0.95,          nan)
(0.6, 0.95,          nan)
(0.7, 0.95,          22 )
(0.8, 0.95,          33 )
  };

\addplot3[only marks,scatter,mark=o,dashed,mark size=5] coordinates{
(0.5, 0.8 ,          33)  
(0.5, 0.9 ,          33)
(0.5, 0.95,          33)
(0.6, 0.95,          33)
};
 \end{axis}
\end{tikzpicture}
\endpgfgraphicnamed
\hfill
\beginpgfgraphicnamed{CVaR_ML_Rev4-fig-beta2}
\begin{tikzpicture}
\begin{axis}[%
 width=0.35\linewidth,
 height=0.30\linewidth,
 xlabel=$x$,
 title=$u$,
 legend style={at={(0.5,0.97)},anchor=north},
 xmin=0.4,xmax=0.6,
 ymin=48,ymax=100,
]
\addplot+[no marks,] table[header=true,x=x,y=u05]{u_beta.dat}; \addlegendentry{$\beta=0.5$};
\addplot+[no marks,] table[header=true,x=x,y=u08]{u_beta.dat}; \addlegendentry{$\beta=0.8$};
\addplot+[no marks,] table[header=true,x=x,y=u09]{u_beta.dat}; \addlegendentry{$\beta=0.9$};
\addplot+[no marks,green!50!black] table[header=true,x=x,y=u095]{u_beta.dat}; \addlegendentry{$\beta=0.95$};
\end{axis}
\end{tikzpicture}
\endpgfgraphicnamed
\caption{Left: number of Newton iterations for different quantiles $\beta$ and rate of decrease of $\varepsilon$ denoted
by $\mu_\varepsilon$. Right: control signals for different $\beta$. Spatial discretization $n_y=129$, other parameters are chosen to equilibrate their individual error contributions. \label{fig:beta}}
\end{figure}

In the last experiment with the 1D PDE, we vary the quantile of the confidence interval in CVaR.
The results are reported in Figure~\ref{fig:beta}.
As $\beta$ increases, the model minimizes the cost $g$ with higher confidence,
which requires a stronger control (see the right plot of Figure~\ref{fig:beta}).
However, the increasing $1/(1-\beta)$ term in the gradient and Hessian renders the Newton method slower and less reliable.
Recall that we start with a large $\varepsilon$ and decrease it geometrically with Newton iterations according to \eqref{eq:mu_iter}.
This allows us to avoid getting too small values of $g''_{\varepsilon}$.
The default value used in the previous experiments was $\mu_{\varepsilon}=0.5$, which was a reasonable balance between the stability of the method and its convergence speed.
However, for larger $\beta$ the Newton method may break at exact machine zeros in $g''_{\varepsilon}$, leading to infinities in the solution increment.
In the left plot of Figure~\ref{fig:beta} we vary both $\beta$ and $\mu_{\varepsilon}$ independently.
Dashed circles denote the experiments where the Newton method failed to converge.
We see that larger $\beta$ requires slower iterations with larger $\mu_{\varepsilon}$.
Of course, too large $\mu_{\varepsilon}$ will just lead to unnecessarily many iterations.
This means that the parameter $\mu_{\varepsilon}$ may need problem-dependent tuning,
although we believe that the observations in Figure~\ref{fig:beta} should serve as a good initial guess in a variety of cases.

%
%
%

\subsection{Elliptic PDE with log-normal coefficient}\label{s:2d}

\begin{figure}[htb]
\centerline{\includegraphics[width=0.4\linewidth]{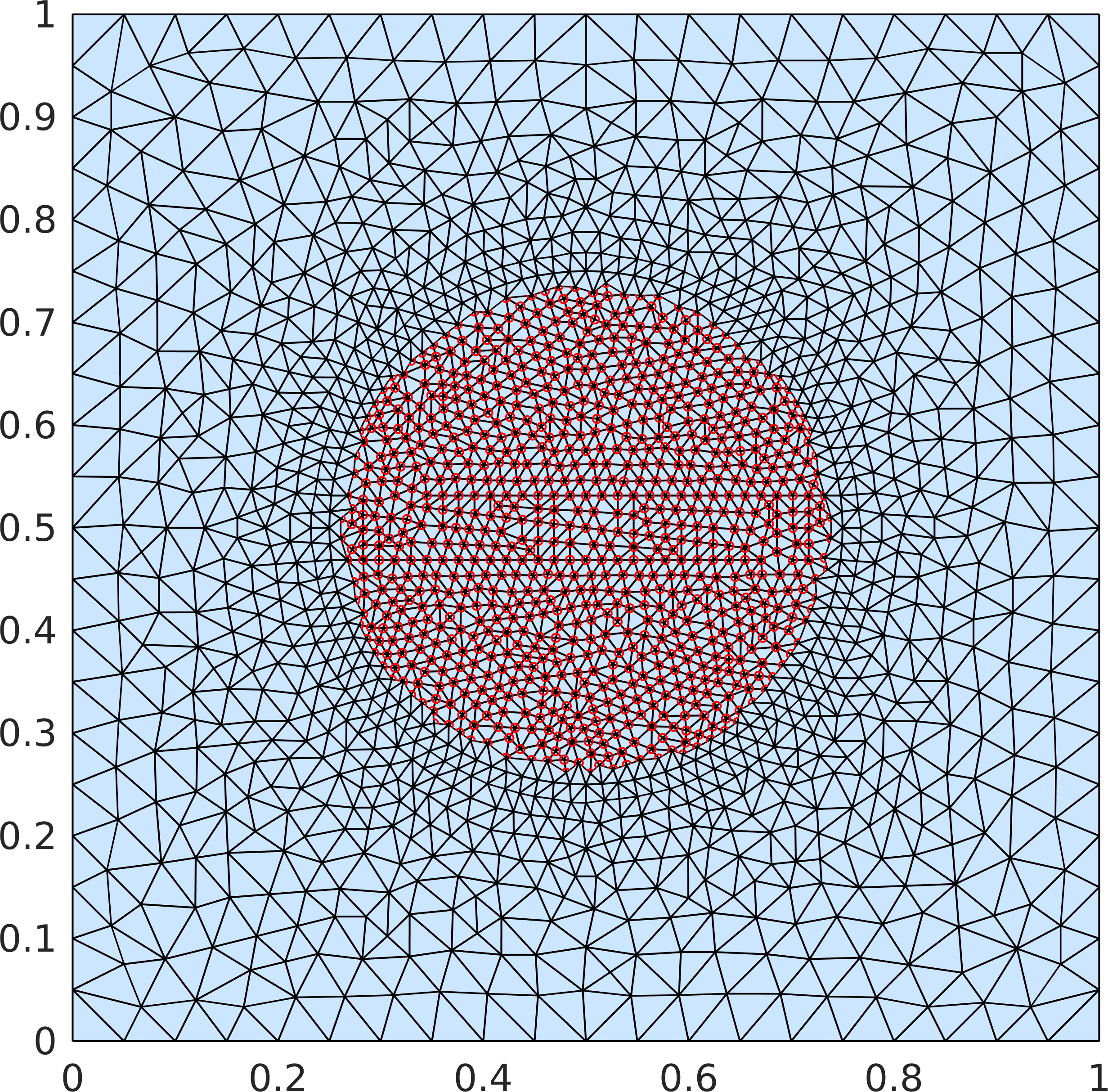}}
\caption{Mesh used to discretize Eq. \eqref{eq:PDE2}. Red points denote grid nodes belonging to the control. \label{fig:2dmesh}}
\end{figure}

Now we consider a larger problem: we solve a PDE on a two-dimensional space,
\begin{align}\label{eq:PDE2}
-\nabla_x \cdot \kappa(x,\xi) \nabla_x y & = Bu, & x \in (0,1)^2, \quad \xi \sim \mathcal{N}(0,\sigma^2 I),\\
y|_{\partial (0,1)^2}& =0, \nonumber
\end{align}
where the coefficient is a log-normal random field.
Namely, we assume that
$$
\log\kappa(x,\xi) = \sum_{k=1}^{d} \sqrt{\lambda_k} \kappa_k(x) \cdot \xi^{(k)},
$$
where $\xi_k \sim \mathcal{N}(0,1)$, the standard normally distributed random variable with $\rho^{(k)}(\xi^{(k)}) = \exp(-(\xi^{(k)})^2/2)/\sqrt{2\pi}$,
is a KL expansion with mean zero and coefficients $\lambda_k,\kappa_k$ defined by the eigenvalue decomposition~\eqref{eq:kleeig} defined by the covariance function
$$
C(x,x') = \sigma^2 \exp\left(-\frac{\|x-x'\|^2_2}{\ell^2}\right), \qquad \ell=0.5, \quad \sigma^2=0.05,
$$
and parametrize $\log\kappa(x,\xi)$ with $d=10$ terms of the KL expansion.
This accounts for $98\%$ of the variance.
The objective function is similar to~\eqref{eq:J1d},
\begin{equation}\label{eq:J2d}
 \mathcal{J}(y,u;\xi) = \frac{1}{2}\|y(u; x,\xi) - y_d(x)\|_{L_2((0,1)^2)}^2, \qquad y_d(x) \equiv 1.
\end{equation}
The control $u$ is defined on a disk inside the domain, and $B$ is the identity insertion operator:
$$
Bu(x) = \left\{\begin{array}{ll}u(x), & \|x-[0.5,0.5]\|_2 \le 0.25, \\ 0, & \mbox{otherwise.}\end{array}\right.
$$
This gives a challenging enough case of incomplete control.
The PDE is discretised using piecewise linear finite elements on a triangular grid with $1829$ nodes, which is shown in Figure~\ref{fig:2dmesh}.
The random variables $\xi^{(1)},\ldots,\xi^{(d)}$ are discretised by collocation at $n_{\xi}=9$ Gauss-Hermite points.


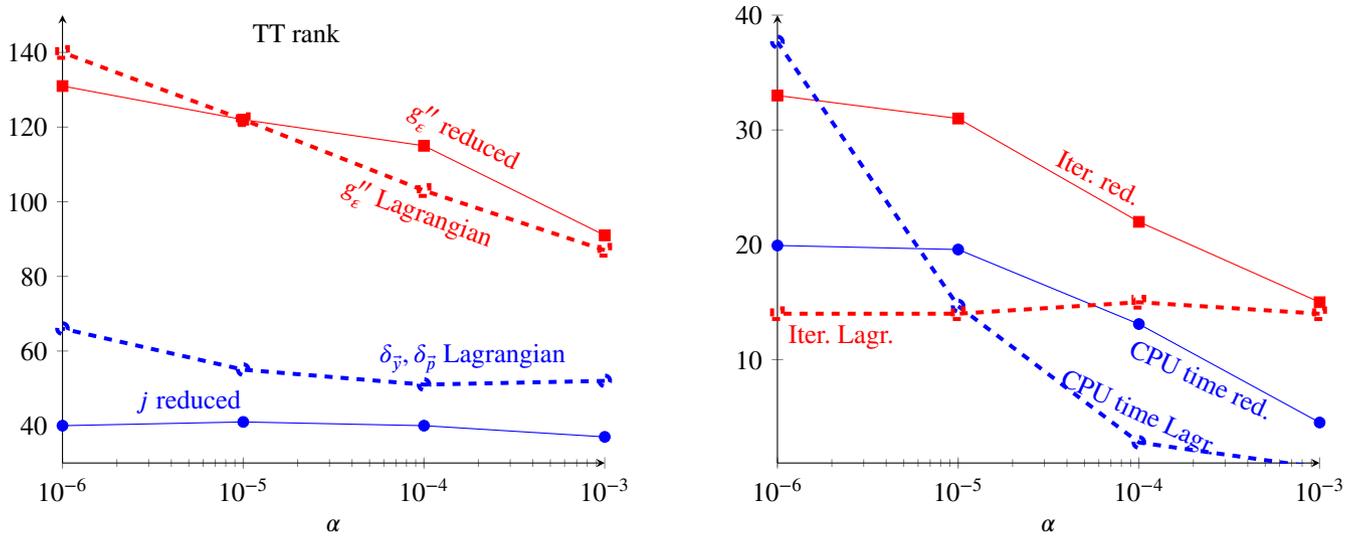
\begin{figure}[tb]
\beginpgfgraphicnamed{CVaR_ML_Rev4-fig-2d1}
\begin{tikzpicture}
\begin{axis}[%
 width=0.40\linewidth,
 height=0.33\linewidth,
 xlabel=$\alpha$,
 xmode=log,
 ymode=normal,
 legend style={at={(0.97,0.03)},anchor=south east},
 ymin=30,ymax=150,
 ]
 \addplot+[blue,solid,mark=*,mark options={style=blue}] coordinates{
 (1e-3, 37)
 (1e-4, 40)
 (1e-5, 41)
 (1e-6, 40)
 } node [pos=0.8,above] {$j$ reduced};
 \addplot+[red,solid,mark=square*,mark options={style=red}] coordinates{
 (1e-3, 91 )
 (1e-4, 115)
 (1e-5, 122)
 (1e-6, 131)
 } node [pos=0.5,above,rotate=-27] {$g_{\varepsilon}''$ reduced};
 \addplot+[blue,dashed,mark=o,mark options={style=blue},line width=1.5pt] coordinates{
 (1e-3, 52)
 (1e-4, 51)
 (1e-5, 55)
 (1e-6, 66)
 } node [pos=0.1,above] {$\delta_{\vec{y}},\delta_{\vec{p}}$ Lagrangian};
 \addplot+[red,dashed,mark=square,mark options={style=red},line width=1.5pt] coordinates{
 (1e-3, 87 )
 (1e-4, 103)
 (1e-5, 122)
 (1e-6, 140)
 } node [pos=0.3,below,rotate=-20] {$g_{\varepsilon}''$ Lagrangian};
 \node[anchor=north west] at (axis cs: 1e-5,150) {TT rank};
\end{axis}
\end{tikzpicture}
\endpgfgraphicnamed
\hfill
\beginpgfgraphicnamed{CVaR_ML_Rev4-fig-2d2}
\begin{tikzpicture}
\begin{axis}[%
 width=0.40\linewidth,
 height=0.33\linewidth,
 xlabel=$\alpha$,
 xmode=log,
 ymode=normal,
 legend style={at={(0.97,0.03)},anchor=south east},
 ymin=1,ymax=40,
 ]
 \addplot+[blue,solid,mark=*,mark options={style=blue}] coordinates{
 (1e-3, 16292/3600)
 (1e-4, 47160/3600)
 (1e-5, 70521/3600)
 (1e-6, 71782/3600)
 } node [pos=0.3,below,rotate=-25] {CPU time red.};
 \addplot+[red,solid,mark=square*,mark options={style=red}] coordinates{
 (1e-3, 15)
 (1e-4, 22)
 (1e-5, 31)
 (1e-6, 33)
 } node [pos=0.5,above,rotate=-27] {Iter. red.};
 \addplot+[blue,dashed,mark=o,mark options={style=blue},line width=1.5pt] coordinates{
 (1e-3, 2110  /3600)
 (1e-4, 10029 /3600)
 (1e-5, 52780 /3600)
 (1e-6, 135577/3600)
 } node [pos=0.1,above,rotate=-25] {CPU time Lagr.};
 \addplot+[red,dashed,mark=square,mark options={style=red},line width=1.5pt] coordinates{
 (1e-3, 14)
 (1e-4, 15)
 (1e-5, 14)
 (1e-6, 14)
 } node [pos=1.0,below,rotate=0,anchor=north west] {Iter. Lagr.};
\end{axis}
\end{tikzpicture}
\endpgfgraphicnamed
\caption{TT ranks (left), CPU times (hours) and numbers of Newton iterations (right) in the reduced and Lagrangian TTRISK formulations depending on the control regularization parameter $\alpha$ (as defined in \eqref{modprob_intro}) in the 2D PDE problem. \label{fig:2d}}
\end{figure}

In this test we set $\beta=0.8$, and correspondingly $\mu_{\varepsilon}=0.7$ (cf.~Figure~\ref{fig:beta}).
The final smoothness width $\varepsilon=3\cdot 10^{-3}$, consistent with
the TT approximation error threshold $\mbox{tol}=10^{-3}$,
the discretization, and KL truncation errors.
In Figure~\ref{fig:2d} we vary the control regularization parameter $\alpha$,
and compare the reduced formulation introduced in section~\ref{sec:reduced},
and the Lagrangian formulation from section~\ref{sec:lagrangian}.

We see that the TT ranks of derivatives of the sigmoid function $g_{\varepsilon}$ in both formulations are of the same scale.
However, since the forward model involves solving a PDE,
the bottleneck is actually the computation of a surrogate of the PDE solution: the TT approximation cost function $j(u;\xi)$ in the reduced formulation, and the block TT format of $\delta_{\vec{y}}$ and $\delta_{\vec{p}}$ in the Lagrangian formulation.
We see that the TT ranks of the Lagrangian solution are 50\% larger than those of $j$, which is actually lower than a factor of $2$ expected for a TT format representing two components simultaneously ($\delta_{\vec{y}}$ and $\delta_{\vec{p}}$).
We see also that the TT ranks of $j$, $\delta_{\vec{y}}$ and $\delta_{\vec{p}}$ depend little on $\alpha$.

Similarly, the number of Newton iterations in the Lagrangian formulation stays nearly constant, whereas in the reduced formulation the number of iterations grows with $\log \alpha$.
This is due to the KKT matrix being the exact Hessian of the Lagrangian,
whereas the reduced formulation uses only an approximate (fixed point) Hessian of $j(u,\xi)$.
This is even more crucial for a large $\beta$ as we have also seen in the previous test.
The imperfection of the reduced Hessian can be seen in the growing number of Newton iterations.
However, if we look at CPU times (Fig.~\ref{fig:2d} right), we notice that the time of the Lagrangian method grows rapidly with $\alpha$ despite nearly constant number of Newton iterations and TT ranks.
This is due to the growth of GMRES iterations in solving Eq.~\eqref{eq:GN}.
The condition number of the KKT matrix preconditioned with Eq.~\eqref{eq:prec} may still deteriorate with $\alpha$ in the considered case of a control on a subdomain.
A better preconditioner may reduce the complexity.
Nevertheless, for moderate values of $\alpha$ even with the current preconditioner the Lagrangian formulation is preferable to the reduced one.

\subsection{Infection ODE model}\label{s:infection}
In the last example, we experiment with an epidemiological ODE model from \cite{DGKP-SEIR-2021},
which was used to estimate the progression of COVID-19 in the UK for $90$ days starting from the 1st March 2020.
The model considers population dynamics split into the following compartments:
\begin{itemize}
 \item Susceptible ($S$): individuals who are not in contact with the virus at the moment.
 \item Exposed ($E$) to the virus, but not yet infectious.
 \item Infected SubClinical ($I^{SC1}$) at the moment, but who may require hospitalizations.
 \item Infected SubClinical ($I^{SC2}$), but recovering without any medical intervention.
 \item Infected Clinical ($I^{C1}$), individuals in the hospital, who may decease after some time.
 \item Infected Clinical ($I^{C2}$) in the hospital, but recovering.
 \item Recovered ($R$) and immune to reinfections.
 \item Deceased ($D$).
\end{itemize}
Each of those categories is further stratified into 5 age groups: 0-19, 20-39, 40-59, 60-79 and 80+.
The age group is denoted by a subscript, e.g. $E_i$ is the number of exposed individuals in the $i$th age group ($i=1,\ldots,5$), and similarly for other compartments.
Each variable of the list above, such as $E=(E_1,\ldots, E_5)$, is thus a vector of length 5.

On the other hand, three simplifications are in order.
First, in the early stage of the pandemic the number of individuals affected by the virus is a small fraction of the entire population.
This allows us to take $S$ constant equal to the initial population, and exclude it from the ODE system.
Similarly, $R$ and $D$ are terminal states in the sense that none of the other variables depend on them.
This again allows us to decouple $R$ and $D$ from the system of ODEs.
Thus, we arrive at a linear system involving only Exposed and Infected numbers:
\begin{align}\label{eq:SIR}
\frac{d}{d\tau}\begin{bmatrix}E \\ I^{SC1} \\ I^{SC2} \\ I^{C1} \\ I^{C2}\end{bmatrix}
& - \begin{bmatrix}
     -\kappa \mathtt{I} & A_u & A_u & 0 & 0 \\
     \kappa \cdot \mathrm{diag}(\rho) & -\gamma_C \mathtt{I} & 0 & 0 & 0 \\
     \kappa \cdot \mathrm{diag}(1-\rho) & 0 & -\gamma_R \mathtt{I} & 0 & 0 \\
     0 & \gamma_C \cdot \mathrm{diag}(\rho') & 0 & -\nu \mathtt{I} & 0 \\
     0 & \gamma_C \cdot \mathrm{diag}(1-\rho') & 0 & 0 & -\gamma_{R,C} \mathtt{I}
    \end{bmatrix}
    \begin{bmatrix}E \\ I^{SC1} \\ I^{SC2} \\ I^{C1} \\ I^{C2}\end{bmatrix} = 0.
\end{align}
Here $\mathtt{I}$ is the identity matrix of size $5$, $\mathrm{diag}(\cdot)$ stretches a vector into a diagonal matrix,
$A_u = \chi \cdot \mathrm{diag}(S) \cdot C_u \cdot \mathrm{diag}(\frac{1}{N})$,
and the remaining variables are model parameters:
\begin{itemize}
 \item $\chi$: probability of contact between $S$ and $I^{SC}$ individuals.
 \item $\kappa=1/d_L$: transition rate of Exposed becoming SubClinical. $d_L$ is the number of days in the Exposed state.
 \item $\gamma_C = 1/d_C$: transition rate of SubClinical turning into Clinical. Similarly, $d_C$ is the number of days this transition takes.
 \item $\gamma_R = 1/d_R$: recovery rate from $I^{SC2}$.
 \item $\gamma_{R,C} = 1/d_{R,C}$: recovery rate from $I^{C2}$.
 \item $\nu = 1/d_D$: death rate from $I^{C1}$.
 \item $\rho=(\rho_1,\ldots,\rho_5)^\top \in \mathbb{R}^5$: age-dependent probabilities of Exposed turning into the first SubClinical category.
 \item $\rho'=(\rho'_1,\ldots,\rho'_5)^\top \in \mathbb{R}^5$: age-dependent probabilities of SubClinical becoming the first Clinical category.
 \item $N = (N_1,\ldots,N_5)^\top \in \mathbb{R}^5$: population sizes in each age group.
 \item $C_u \in \mathbb{R}^{5\times 5}$: the contact matrix, which depends on the control $u$.
\end{itemize}
Scalar-vector operations ($1/N$, $1-\rho$, etc.) are understood as elementwise operations.

Another parameter is the number of infected individuals on day 0 (1st March) $N^{0}$.
It is split further across the age groups as follows:
$$
N^{in}:=(N^{in}_1, N^{in}_2, N^{in}_3, N^{in}_4, N^{in}_5)^\top = (0.1, 0.4, 0.35, 0.1, 0.05)^\top N^{0}.
$$
The ODE \eqref{eq:SIR} is initialized by setting
$$
E(0) = \frac{N^{in}}{3}, \quad I^{SC1}(0) = \frac{2}{3} \mathrm{diag}(\rho) N^{in}, \quad I^{SC2}(0) = \frac{2}{3} \mathrm{diag}(1-\rho) N^{in}, \quad I^{C1}(0) = I^{C2}(0)=0.
$$
The population size $S=N$ is taken from the office of national statistics, mid 2018 estimate.

The contact matrix consists of four contributions, corresponding to people activities\footnote{Note that in this section $\alpha$ is the notation from the original paper \cite{DGKP-SEIR-2021}, not a regularization parameter.}:
\begin{equation}
 C_u = \mathrm{diag}(\alpha^{home}) C^{home} + \mathrm{diag}(\alpha^{work}) C^{work} + \mathrm{diag}(\alpha^{school}) C^{school} + \mathrm{diag}(\alpha^{other}) C^{other},
\end{equation}
where $C^{*}$ are pre-pandemic contact matrices (for details of their estimation see \cite{DGKP-SEIR-2021}),
and $\alpha^{*}$ are coefficients of reduction of activities introduced in March 2020.
(Here $*$ stands for $home,work,school$ or $other$.)
In turn, those are constructed as follows.
Firstly, we set $\alpha^{home} = (1,\ldots,1)$, since home contacts cannot be influenced.
For the remaining activities, noting that the lockdown in the UK was called on day 17 (18th March), we set
\begin{equation}
\alpha^*(\tau) = \left\{\begin{array}{ll}(1,1,1,1,1)^\top, & \tau<17, \\ (\alpha_{123}(1-u^*(\tau)) ,\alpha_{123} (1-u^*(\tau)),\alpha_{123} (1-u^*(\tau)),\alpha_4,\alpha_5)^\top, & \mbox{otherwise}, \end{array}\right.
\end{equation}
where $u^{work}$, $u^{school}$ and $u^{other}$ are the control signals (the lockdown measures) that we are going to optimize,
and $\alpha_{123},\alpha_4,\alpha_5$ are their proportions in the corresponding age groups.

In the cost function we penalize the total fatalities and the hospital capacity exceedance.
As long as \eqref{eq:SIR} is solved,
the number of deaths can be calculated directly as
\begin{equation}\label{eq:deaths}
D(\tau) = \nu \int_{0}^{\tau} I^{C1}(s) ds.
\end{equation}
Moreover, we aggregate $I^{C} = \sum_{i=1}^5 I^{C1}_i + I^{C2}_i$ and penalize $I^C$ exceeding $10000$.
Finally, we penalize the strength of the lockdown measures, i.e. the norm of the control $u(\tau) = (u^{work}(\tau), u^{school}(\tau), u^{other}(\tau))$,
as well as constraining $u^{work} \in [0, 0.69]$, $u^{school} \in [0, 0.9]$ and $u^{other} \in [0, 0.59]$.
This gives us a deterministic cost
\begin{equation}\label{eq:cost-det}
 \mathcal{J}^{Det}(u) = \frac{1}{2}\left[ D(T) + \int_0^T \max(I^C(\tau) - 10000, 0)d\tau + \zeta \int_{17}^{T} \|u(\tau)\|_2^2 d\tau \right],
\end{equation}
where $T=90$ is the simulation interval, and $\zeta$ is the regularization parameter.

However, optimization of \eqref{eq:cost-det} may be misleading, since the model parameters are not known in advance, and can only be estimated.
In particular, \cite{DGKP-SEIR-2021} employed an Approximate Bayesian Computation (ABC), which used existing observations of daily deaths and hospitalizations in the UK to form the likelihood, and consequently the posterior probability density function.
This renders model parameters into random variables, which are distributed according to the posterior density.
In turn, this motivates a modification of \eqref{eq:cost-det} into a risk-averse cost function, for example, using CVaR with $\beta=0.5$.
More precisely, we aim to minimize
\begin{equation}\label{eq:cost-cvar}
 \mathcal{J}^{CVaR}(u) = \mathrm{CVaR}^{\varepsilon}_{0.5}\left(\frac{1}{2}D(T) + \frac{1}{2}\int_0^T \max(I^C(\tau) - 10000, 0)d\tau \right) + \frac{\zeta}{2} \int_{17}^{T} \|u(\tau)\|_2^2 d\tau.
\end{equation}
Ideally, the expectations in CVaR need to be computed with respect to the posterior density from ABC.
However, the latter is a complicated multivariate function,
which lacks an independent variable parametrization necessary to set up the discretization and TT approximation.
(A possible solution to this using optimal transport \cite{CD-DIRT-2021} can be a matter of future research.)
As a proof of concept, we simplify the distribution to independent uniform, reflecting means and variances estimated by ABC.
Thus, we assume
\begin{align}\label{eq:SIR-params}
\chi          & \sim \mathcal{U}(0.13- 0.03\sigma, 0.13+ 0.03\sigma), & d_L           & \sim \mathcal{U}(1.57- 0.42\sigma, 1.57+ 0.42\sigma), \\\nonumber
d_C           & \sim \mathcal{U}(2.12 -0.80\sigma, 2.12 +0.80\sigma), & d_R           & \sim \mathcal{U}(1.54 -0.40\sigma, 1.54 +0.40\sigma), \\\nonumber
d_{R,C}       & \sim \mathcal{U}(12.08-1.51\sigma, 12.08+1.51\sigma), & d_D           & \sim \mathcal{U}(5.54 -2.19\sigma, 5.54 +2.19\sigma), \\\nonumber
\rho_1        & \sim \mathcal{U}(0.06 -0.03\sigma, 0.06 +0.03\sigma), & \rho_2        & \sim \mathcal{U}(0.05 -0.03\sigma, 0.05 +0.03\sigma), \\\nonumber
\rho_3        & \sim \mathcal{U}(0.08 -0.04\sigma, 0.08 +0.04\sigma), & \rho_4        & \sim \mathcal{U}(0.54 -0.22\sigma, 0.54 +0.22\sigma), \\\nonumber
\rho_5        & \sim \mathcal{U}(0.79 -0.14\sigma, 0.79 +0.14\sigma), & \rho'_1       & \sim \mathcal{U}(0.26 -0.23\sigma, 0.26 +0.23\sigma), \\\nonumber
\rho'_2       & \sim \mathcal{U}(0.28 -0.25\sigma, 0.28 +0.25\sigma), & \rho'_3       & \sim \mathcal{U}(0.33 -0.27\sigma, 0.33 +0.27\sigma), \\\nonumber
\rho'_4       & \sim \mathcal{U}(0.26 -0.11\sigma, 0.26 +0.11\sigma), & \rho'_5       & \sim \mathcal{U}(0.80 -0.13\sigma, 0.80 +0.13\sigma), \\\nonumber
N^{0}         & \sim \mathcal{U}(276  - 133\sigma, 276  +133 \sigma), & \alpha_{123}  & \sim \mathcal{U}(0.63 -0.21\sigma, 0.63 +0.21\sigma), \\\nonumber
\alpha_4      & \sim \mathcal{U}(0.57 -0.23\sigma, 0.57 +0.23\sigma), & \alpha_5      & \sim \mathcal{U}(0.71 -0.23\sigma, 0.71 +0.23\sigma),   \nonumber
\end{align}
where $\sigma$ is a variance tuning parameter.
That is, \eqref{eq:SIR-params} form a random vector
$$
\xi=(\chi,d_L,d_C,d_R,d_{R,C},d_D,\rho_1,\rho_2,\rho_3\rho_4,\rho_5,\rho'_1,\rho'_2,\rho'_3\rho'_4,\rho'_5,N^{0},\alpha_{123},\alpha_4,\alpha_5)
$$
of dimension $d=20$,
the state vector is
$$
y=(E_1,\ldots,E_5,~I^{SC1}_1,\ldots,I^{SC1}_5,~I^{SC2}_1,\ldots,I^{SC2}_5,~I^{C1}_1,\ldots,I^{C1}_5,~I^{C2}_1,\ldots,I^{C2}_5),
$$
and
the left hand side of~\eqref{eq:SIR} constitutes the constraints $c$.

For the deterministic optimization \eqref{eq:cost-det} we set all variables to their means.
For the stochastic optimization \eqref{eq:cost-cvar} we discretize each variable by a $3$-point Gauss-Legendre quadrature rule on the corresponding interval,
and the TT approximations are computed with relative error threshold of $10^{-2}$.
The control is discretized with a Gauss-Lagrange interpolation on $[17,T]$ with $7$ points, which makes the total dimension of the discretized control space $21$.
This allows us to use the reduced optimization formulation.
Moreover, since the cost contains non-smooth functions, we abandon the Newton scheme,
resorting to the Projected Gradient Descent method with a finite difference computation of the gradient of \eqref{eq:cost-cvar}.
The ODE \eqref{eq:SIR} is solved using the implicit Euler method with time step $0.1$.
The control regularization $\zeta=100$ is taken from \cite{DGKP-SEIR-2021},
and the CVaR smoothing parameter $\varepsilon=1000$ corresponds to the relative width of the smoothed region $\varepsilon/t < 10^{-2}$, matching the bias of the smoothed CVaR and the TT approximation error.

\begin{table}[t]
\centering
\caption{CVaR optimizer behavior for different variability of random variables \label{tab:SIR}}
\begin{tabular}{c|ccc|cc}
 $\sigma$ & $\mathcal{J}^{CVaR}$ & $\frac{\zeta}{2} \int_{17}^{T} \|u(\tau)\|_2^2 d\tau$ & $t$ & Iterations & TT rank($\grad{\mathcal{J}^{CVaR}}{u}$) \\ \hline
$10^{-2}$ & 171168 &         3107 & 158382 & 16 & 13 \\
$10^{-1}$ & 289160 &         3439 & 157750 & 20 & 68 \\
\end{tabular}
\end{table}

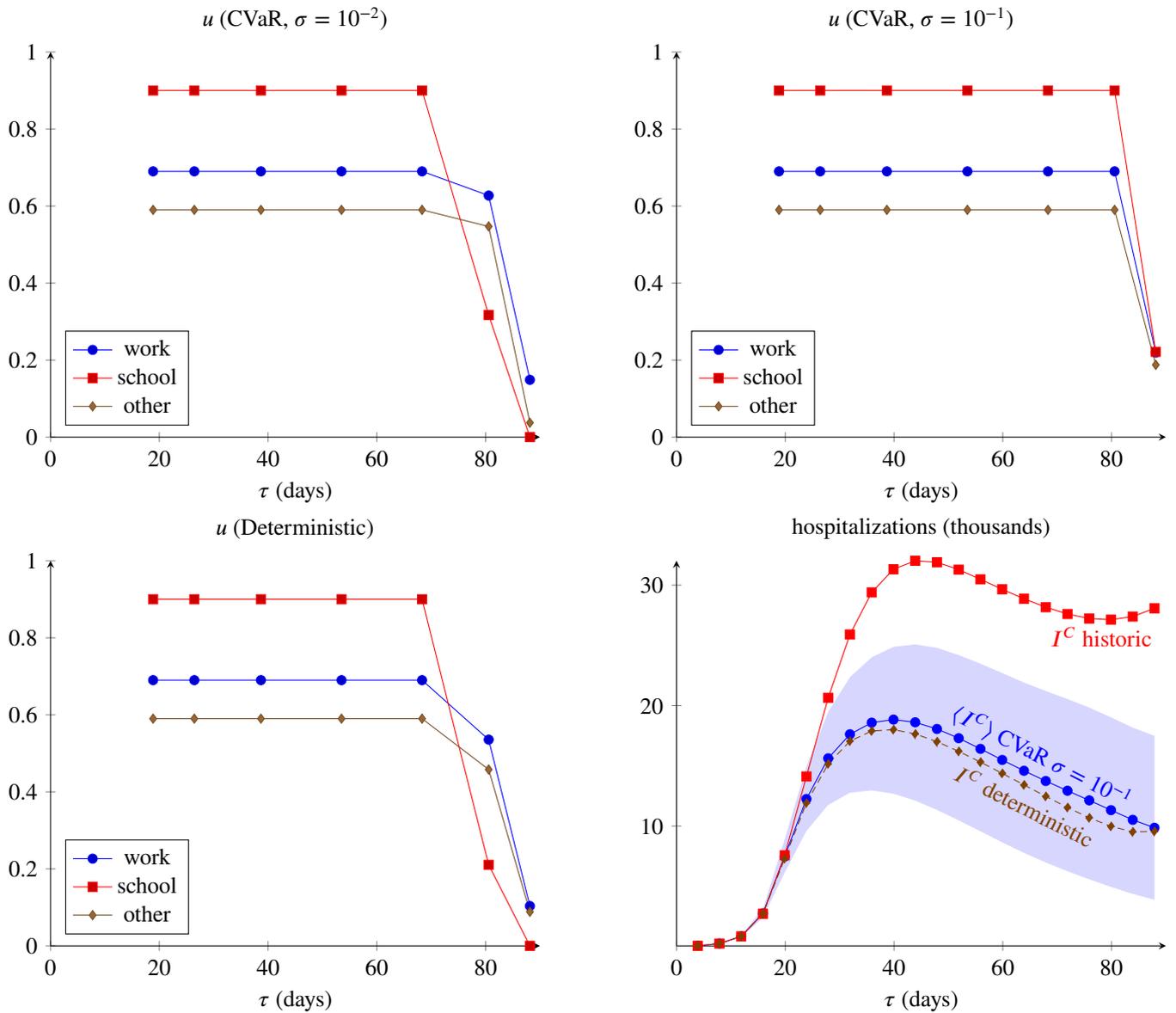
\begin{figure}[htb]
\beginpgfgraphicnamed{CVaR_ML_Rev4-fig-SIR1}
\begin{tikzpicture}
 \begin{axis}[%
 width=0.42\linewidth,
 height=0.33\linewidth,
 xlabel=$\tau$ (days),
 title={$u$ (CVaR, $\sigma=10^{-2}$)},
 xmin=0,xmax=90,
 ymin=0,ymax=1,
 legend style={at={(0.03,0.03)},anchor=south west},
 ]
 \addplot+[] coordinates{
   (18.8576,     0.6900)
   (26.4341,     0.6900)
   (38.6867,     0.6900)
   (53.5000,     0.6900)
   (68.3133,     0.6900)
   (80.5659,     0.6272)
   (88.1424,     0.1488)
 }; \addlegendentry{work};
 \addplot+[] coordinates{
   (18.8576,    0.9000)
   (26.4341,    0.9000)
   (38.6867,    0.9000)
   (53.5000,    0.9000)
   (68.3133,    0.9000)
   (80.5659,    0.3169)
   (88.1424,         0)
 }; \addlegendentry{school};
 \addplot+[mark=diamond*] coordinates{
   (18.8576,    0.5900  )
   (26.4341,    0.5900  )
   (38.6867,    0.5900  )
   (53.5000,    0.5900  )
   (68.3133,    0.5900  )
   (80.5659,    0.5470  )
   (88.1424,    0.0373  )
 }; \addlegendentry{other};
 \end{axis}
\end{tikzpicture}
\endpgfgraphicnamed
\hfill
\beginpgfgraphicnamed{CVaR_ML_Rev4-fig-SIR2}
\begin{tikzpicture}
 \begin{axis}[%
 width=0.42\linewidth,
 height=0.33\linewidth,
 xlabel=$\tau$ (days),
 title={$u$ (CVaR, $\sigma=10^{-1}$)},
 xmin=0,xmax=90,
 ymin=0,ymax=1,
 legend style={at={(0.03,0.03)},anchor=south west},
 ]
 \addplot+[] coordinates{
   (18.8576,    0.6900)
   (26.4341,    0.6900)
   (38.6867,    0.6900)
   (53.5000,    0.6900)
   (68.3133,    0.6900)
   (80.5659,    0.6900)
   (88.1424,    0.2194)
 }; \addlegendentry{work};
 \addplot+[] coordinates{
   (18.8576,    0.9000)
   (26.4341,    0.9000)
   (38.6867,    0.9000)
   (53.5000,    0.9000)
   (68.3133,    0.9000)
   (80.5659,    0.9000)
   (88.1424,    0.2215)
 }; \addlegendentry{school};
 \addplot+[mark=diamond*] coordinates{
   (18.8576,    0.5900)
   (26.4341,    0.5900)
   (38.6867,    0.5900)
   (53.5000,    0.5900)
   (68.3133,    0.5900)
   (80.5659,    0.5900)
   (88.1424,    0.1876)
 }; \addlegendentry{other};
 \end{axis}
\end{tikzpicture}
\endpgfgraphicnamed

\beginpgfgraphicnamed{CVaR_ML_Rev4-fig-SIR3}
\begin{tikzpicture}
 \begin{axis}[%
 width=0.42\linewidth,
 height=0.33\linewidth,
 xlabel=$\tau$ (days),
 title={$u$ (Deterministic)},
 xmin=0,xmax=90,
 ymin=0,ymax=1,
 legend style={at={(0.03,0.03)},anchor=south west},
 ]
 \addplot+[] coordinates{
   (18.8576,    0.6900)
   (26.4341,    0.6900)
   (38.6867,    0.6900)
   (53.5000,    0.6900)
   (68.3133,    0.6900)
   (80.5659,    0.5355)
   (88.1424,    0.1035)
 }; \addlegendentry{work};
 \addplot+[] coordinates{
   (18.8576,    0.9000)
   (26.4341,    0.9000)
   (38.6867,    0.9000)
   (53.5000,    0.9000)
   (68.3133,    0.9000)
   (80.5659,    0.2105)
   (88.1424,         0)
 }; \addlegendentry{school};
 \addplot+[mark=diamond*] coordinates{
   (18.8576,    0.5900)
   (26.4341,    0.5900)
   (38.6867,    0.5900)
   (53.5000,    0.5900)
   (68.3133,    0.5900)
   (80.5659,    0.4579)
   (88.1424,    0.0885)
 }; \addlegendentry{other};
 \end{axis}
\end{tikzpicture}
\endpgfgraphicnamed
\hfill
\beginpgfgraphicnamed{CVaR_ML_Rev4-fig-SIR4}
\begin{tikzpicture}
 \begin{axis}[%
 width=0.42\linewidth,
 height=0.33\linewidth,
 xlabel=$\tau$ (days),
 title={hospitalizations (thousands)},
 xmin=0,xmax=90,
 ]
 \addplot+[no marks,color=white,name path=minus] coordinates{
   ( 3.9000,    0.0615)
   ( 7.9000,    0.2374)
   (11.9000,    0.7623)
   (15.9000,    2.3394)
   (19.9000,    6.0513)
   (23.9000,    9.5330)
   (27.9000,   11.7112)
   (31.9000,   12.7305)
   (35.9000,   12.9277)
   (39.9000,   12.6524)
   (43.9000,   12.0760)
   (47.9000,   11.3054)
   (51.9000,   10.4335)
   (55.9000,    9.5307)
   (59.9000,    8.6058)
   (63.9000,    7.7364)
   (67.9000,    6.9411)
   (71.9000,    6.2150)
   (75.9000,    5.5434)
   (79.9000,    4.9155)
   (83.9000,    4.3379)
   (87.9000,    3.8508)
};
 \addplot+[no marks,color=white,name path=plus] coordinates{
   ( 3.9000,    0.0688)
   ( 7.9000,    0.2754)
   (11.9000,    0.9584)
   (15.9000,    3.2402)
   (19.9000,    9.0312)
   (23.9000,   15.0767)
   (27.9000,   19.5558)
   (31.9000,   22.4201)
   (35.9000,   24.0503)
   (39.9000,   24.9170)
   (43.9000,   25.1119)
   (47.9000,   24.8297)
   (51.9000,   24.2309)
   (55.9000,   23.5187)
   (59.9000,   22.7419)
   (63.9000,   21.9515)
   (67.9000,   21.2491)
   (71.9000,   20.5857)
   (75.9000,   19.8743)
   (79.9000,   19.0689)
   (83.9000,   18.2106)
   (87.9000,   17.5293)
};

\addplot[blue!16!white] fill between[of = minus and plus];

 \addplot+[blue,mark=*] coordinates{
   ( 3.9000,    0.0650)
   ( 7.9000,    0.2564)
   (11.9000,    0.8568)
   (15.9000,    2.7678)
   (19.9000,    7.4794)
   (23.9000,   12.2512)
   (27.9000,   15.6264)
   (31.9000,   17.6198)
   (35.9000,   18.5825)
   (39.9000,   18.8382)
   (43.9000,   18.6130)
   (47.9000,   18.0596)
   (51.9000,   17.2925)
   (55.9000,   16.4086)
   (59.9000,   15.4890)
   (63.9000,   14.5918)
   (67.9000,   13.7433)
   (71.9000,   12.9355)
   (75.9000,   12.1376)
   (79.9000,   11.3245)
   (83.9000,   10.5209)
   (87.9000,    9.8634)
} node [pos=0.75,above,rotate=-27] {$\langle I^C\rangle $ CVaR $\sigma=10^{-1}$};
 \addplot+[red,mark=square*,mark options={style=red}] coordinates{
   ( 3.9000,     0.0649)
   ( 7.9000,     0.2551)
   (11.9000,     0.8480)
   (15.9000,     2.7222)
   (19.9000,     7.5793)
   (23.9000,    14.1194)
   (27.9000,    20.6464)
   (31.9000,    25.8989)
   (35.9000,    29.3929)
   (39.9000,    31.3104)
   (43.9000,    32.0187)
   (47.9000,    31.8938)
   (51.9000,    31.2797)
   (55.9000,    30.4762)
   (59.9000,    29.6473)
   (63.9000,    28.8639)
   (67.9000,    28.1656)
   (71.9000,    27.6001)
   (75.9000,    27.2263)
   (79.9000,    27.1277)
   (83.9000,    27.3871)
   (87.9000,    28.0681)
} node[pos=0.9,below]{$I^C$ historic};
 \addplot+[orange!50!black,mark=diamond*,mark options={style=orange!50!black}] coordinates{
   ( 3.9000,    0.0649)
   ( 7.9000,    0.2551)
   (11.9000,    0.8480)
   (15.9000,    2.7222)
   (19.9000,    7.3049)
   (23.9000,   11.9120)
   (27.9000,   15.1612)
   (31.9000,   17.0495)
   (35.9000,   17.8906)
   (39.9000,   18.0054)
   (43.9000,   17.6495)
   (47.9000,   17.0073)
   (51.9000,   16.2023)
   (55.9000,   15.3107)
   (59.9000,   14.3757)
   (63.9000,   13.4219)
   (67.9000,   12.4694)
   (71.9000,   11.5447)
   (75.9000,   10.6911)
   (79.9000,    9.9783)
   (83.9000,    9.5245)
   (87.9000,    9.5606)
} node [pos=0.75,below,rotate=-27] {$I^C$ deterministic};
 \end{axis}
\end{tikzpicture}
\endpgfgraphicnamed
\caption{Control signals and predicted hospitalizations numbers for different optimization strategies. In the bottom right plot, blue circles denote mean $I^C$ after CVaR optimization, shaded area denotes 95\% confidence interval. \label{fig:SIR}}
\end{figure}

In Table~\ref{tab:SIR} we vary $\sigma$ and investigate the behavior of the reduced TT formulation for CVaR.
Note that we need more iterations and higher TT ranks for the larger variance.
This is also reflected in a larger total cost, which is dominated by the CVaR term.
In Figure~\ref{fig:SIR}, we compare the controls computed with the two values of $\sigma$, as well as the minimizer of the deterministic cost~\eqref{eq:cost-det}.
We see that a small $\sigma$ yields the solution which is similar to the deterministic one.
However, the risk-averse control for a larger variance of the parameters is more conservative: it tends to be larger, hitting the constraints at almost the entire interval except the final point, where the control stops making any influence on the system.

In the bottom right panel of Figure~\ref{fig:SIR} we plot the predicted hospitalization numbers.
The historic numbers are obtained by simulating the deterministic model (using mean values of the parameters in \eqref{eq:SIR-params}) with the control derived from the smoothed Google daily mobility data\footnote{The historic $I^C$ differs from that in \cite[Figure 6]{DGKP-SEIR-2021}. This is most likely due to imprecise reproduction of the parameters. However, the values agree within a factor of $2$, which indicates that the qualitative behaviour of the model is correct.}.
We see that both optimization techniques allow one to reduce the hospital occupancy.
However, the deterministic approach tends to underestimate $I^C$ compared to the mean risk-averse value.
In addition, the CVaR frameworks enables a rigorous uncertainty quantification.

\section{Conclusion and outlook}
This paper has introduced a new algorithm called TTRISK to solve risk-averse optimization problems constrained by differential equations (PDEs or ODEs). TTRISK can be applied to both the reduced and full-space formulations. The article also introduces a control variate correction to get unbiased estimators. Various strategies to choose the underlying algorithmic parameters have been outlined throughout the paper, especially in the numerics section. This is carried out by carefully taking into account all the approximation errors.

The TT framework offers multiple advantages, for instance, our numerical examples illustrate that it can help overcome the so-called ``curse of dimensionality''. Indeed, the approach introduced here has been successfully applied to a realistic problem, with 20 random variables, to study the propagation of COVID-19 and to devise optimal lockdown strategies.

This article aims to initiate new research directions in the field of risk-averse optimization. There are many open questions that one could pursue, a few examples are:
\begin{enumerate}
\item Convergence analysis of TTRISK in both reduced and full-space formulations.

\item Convergence analysis of TTRISK in the presence of control variate terms.

\item Convergence analysis of preconditioned Gauss-Newton method (cf.~Section~\ref{s:precond})
 for the full-space formulation.
\item The Gauss-Newton system (cf.~Section~\ref{s:GN}) and preconditioning (cf.~Section~\ref{s:precond})
   for the full space formulation eliminates the control and it considers a formulation in terms
   of $y$, $p$, and $t$. It maybe interesting to design preconditioners which can handle
 control $u$ directly in the full-space formulation.
\end{enumerate}

\section*{Acknowledgement}
HA and AO are partially supported by NSF grants DMS-2110263, DMS-1913004 
and the Air Force Office of Scientific Research under Award NO: FA9550-19-1-0036 and FA9550-22-1-0248.
SD is thankful for the support from Engineering and Physical Sciences Research 
Council (EPSRC) New Investigator Award EP/T031255/1 and New Horizons grant EP/V04771X/1.

\noindent
\textbf{Data Access.} Matlab codes for reproducing the experiments are available at \url{https://github.com/dolgov/TTRISK}

\noindent
\textbf{Conflict of interest.} This study does not have any conflicts to disclose.

\bibliography{refs}

\appendix

\section{Cross approximation}\label{sec:cross}
Suppose $f(\xi)$ is defined as a minimizer of a cost function $\mathcal{C}(f)$.
We aim at minimizing $\mathcal{C}$ over a TT decomposition instead by driving $\nabla \mathcal{C}(\tilde f)=0$.
However, $\tilde f$ is a product expansion in the actual degrees of freedom $\{\mathbf{F}^{(k)}\}$, which makes the problem nonlinear even if $\grad{\mathcal{C}(f)}{f}$ was linear.
Instead of running a generic nonlinear optimization of all TT cores simultaneously, we can resort to
the Alternating Least Squares~\cite{holtz-ALS-DMRG-2012} and Density Matrix Renormalization Group~\cite{white-dmrg-1993,schollwock-2005} methods.
Those alleviate the nonlinearity issue by iterating over $k=1,\ldots,d$,
solving only $\grad{\mathcal{C}(\tilde f)}{\mathbf{F}^{(k)}}=0$ in each step, similarly to the coordinate descent method.
Note that $\tilde f$ is linear in each \emph{individual} factor $\mathbf{F}^{(k)}$.

Applying this coordinate descent idea to the problem of interpolating a given function with a TT decomposition yields a family of \emph{TT-Cross} algorithms~\cite{ot-ttcross-2010,so-dmrgi-2011proc,ds-parcross-2020}.
Suppose we are given a \emph{procedure} to evaluate a continuous function $f(\xi)$ at any given $\xi$.
We iterate over $k=1,\ldots,d$ and in each step compute $\mathbf{F}^{(k)}$ by solving an interpolation condition
\begin{equation}\label{eq:ttinterp}
\tilde f(\xi) = f(\xi) \qquad \forall \xi \in \boldsymbol\Xi_k:=\{(\xi^{(1)}_j,\ldots,\xi^{(d)}_j) \in \mathbb{R}^d: j=1,\ldots,r_{k-1}n_{\xi}r_k\},
\end{equation}
over some carefully chosen sampling sets $\boldsymbol\Xi_k$.
Stretching the tensor $\mathbf{F}^{(k)}$ into a vector $f^{(k)} = [\mathbf{F}^{(k)}(s_{k-1},i,s_k)] \in \mathbb{R}^{r_{k-1}n_{\xi} r_k}$,
we can write \eqref{eq:ttinterp} as a linear equation $F_{\neq k} f^{(k)} = f(\boldsymbol\Xi_k)$, where each row of the matrix $F_{\neq k} \in \mathbb{R}^{r_{k-1} n_{\xi} r_k \times r_{k-1} n_{\xi} r_k}$ is given by
\begin{equation}\label{eq:cross-frame}
 F_{\neq k}(j, s_{k-1} i s_k) = \sum_{s_0,\ldots,s_{k-2}} F^{(1)}_{s_0,s_1}(\xi^{(1)}_j) \cdots F^{(k-1)}_{s_{k-2},s_{k-1}}(\xi^{(k-1)}_j) \cdot \ell_i(\xi^{(k)}_j) \cdot \sum_{s_{k+1},\ldots,s_d} F^{(k+1)}_{s_{k},s_{k+1}}(\xi^{(k+1)}_j) \cdots F^{(d)}_{s_{d-1},s_d}(\xi^{(d)}_j),
\end{equation}
where $\xi_j \in \boldsymbol\Xi_k$, $j=1,\ldots,r_{k-1}n_{\xi}r_k$.
To compute this efficiently in a recursive manner similar to \eqref{eq:int-rec}, we restrict $\boldsymbol\Xi_k$ to the Cartesian form
\begin{equation}\label{eq:set-kron}
 \boldsymbol\Xi_{k} = \Xi_{<k} \times \Xi^{(k)} \times \Xi_{>k},
\end{equation}
for some chosen sets
$$
\Xi_{<k} = \{(\xi^{(1)}_{s_{k-1}},\ldots,\xi^{(k-1)}_{s_{k-1}})\}_{s_{k-1}=1}^{r_{k-1}}, \qquad \Xi_{>k} = \{(\xi^{(k+1)}_{s_k},\ldots,\xi^{(d)}_{s_k})\}_{s_k=1}^{r_k},
$$
including $\Xi_{<1} = \Xi_{>d} = \emptyset$,
and $\Xi^{(k)} = \{\xi^{(k)}_i\}$ are the quadrature nodes associated with $\{\ell_i(\xi^{(k)})\}$.
This allows us to write \eqref{eq:cross-frame} in the form
\begin{equation}\label{eq:frame-kron}
F_{\neq k} = F_{<k} \otimes \ell(\Xi^{(k)}) \otimes F_{>k},
\end{equation}
where $\otimes$ denotes the Kronecker product operator and
\begin{align}
F_{<k} & = \left[\sum_{s_0,\ldots,s_{k-2}} F^{(1)}_{s_0,s_1}(\xi^{(1)}_j) \cdots F^{(k-1)}_{s_{k-2},s_{k-1}}(\xi^{(k-1)}_j)\right] \in \mathbb{R}^{r_{k-1} \times r_{k-1}}, & (\xi^{(1)}_j,\ldots,\xi^{(k-1)}_j) & \in \Xi_{<k}, \\
F_{>k} & = \left[\sum_{s_{k+1},\ldots,s_d} F^{(k+1)}_{s_{k},s_{k+1}}(\xi^{(k+1)}_j) \cdots F^{(d)}_{s_{d-1},s_d}(\xi^{(d)}_j)\right] \in \mathbb{R}^{r_k \times r_k}, & (\xi^{(k+1)}_{j},\ldots,\xi^{(d)}_{j}) & \in \Xi_{>k}.
\end{align}
Moreover, for the \emph{left} and \emph{right} sets $\Xi_{<k},\Xi_{>k}$ we assume \emph{nestedness} conditions
\begin{equation}\label{eq:nest}
\Xi_{<k+1} \subset \Xi_{<k} \times \Xi^{(k)}, \qquad \Xi_{>k-1} \subset \Xi^{(k)} \times \Xi_{>k}.
\end{equation}
This way, given $F_{<k}$ or $F_{>k}$, we can continue the recursion by computing
\begin{align}\label{eq:cross-ileft}
 \overline{F}_{\le k}(s_{k-1} i, s_k) & = \sum_{t_{k-1}} F_{<k}(s_{k-1},t_{k-1}) F^{(k)}_{t_{k-1},s_k}(\xi^{(k)}_i),
 \quad \mbox{and} \\
 \overline{F}_{\ge k}(s_{k-1}, i s_k) & = \sum_{t_{k}} F^{(k)}_{s_{k-1},t_k}(\xi^{(k)}_i) F_{<k}(t_k,s_k),\label{eq:cross-iright}
\end{align}
and extracting $F_{<k+1}$ and $F_{>k-1}$ simply as submatrices of $\overline{F}_{\le k}$ and $\overline{F}_{\ge k}$, respectively.
This needs $\mathcal{O}(n_{\xi} r^3)$ operations.

The particular elements extracted from e.g. $\Xi_{<k} \times \Xi^{(k)}$ (and $\overline{F}_{\le k}$) are chosen to ensure that $F_{<k+1}$ (and hence $F_{\neq k+1}$) are as well conditioned as possible.
We use the Maximum Volume (\emph{maxvol}) algorithm \cite{gostz-maxvol-2010} which performs an iterative optimization of the volume (absolute determinant) of the submatrix $F_{<k+1} = \overline{F}_{\le k}(\mathcal{I}_{\le k},:)$ over the index set $\mathcal{I}_{\le k} = \{s_{k-1},i\}$ of cardinality $r_k$.
We can apply the same algorithm to $\overline{F}_{\ge k}^\top$ to find an index set $\mathcal{I}_{\ge k} = \{i,s_k\}$, with cardinality $\#\mathcal{I}_{\ge k} = r_{k-1}$, that provides $|\mathrm{det}\overline{F}_{\ge k}(:,\mathcal{I}_{\ge k})| \approx \max_{\mathcal{I}} |\mathrm{det}\overline{F}_{\ge k}(:,\mathcal{I})|$.
The recursion over the sampling sets \eqref{eq:nest} is arranged by collecting $\Xi_{<k+1} = \{\Xi_{<k}(s_{k-1}), \xi^{(k)}_{i}\}$ for $s_{k-1},i \in \mathcal{I}_{\le k}$, and similarly for $\Xi_{>k-1}$.
The entire iteration, which we call TT-cross, is outlined in Algorithm~\ref{alg:tt-cross}.
\begin{algorithm}[htb]
\caption{TT-Cross}
\label{alg:tt-cross}
\begin{algorithmic}[1]
 \State Choose initial sets $\Xi_{<k}$, $k=2,\ldots,d$, stopping threshold $\delta>0$.
 \While{first iteration or $\|\tilde f(\xi) - \tilde f_{\rm prev}(\xi) \|>\delta \|\tilde f(\xi)\|$}
   \For{$k=d,d-1,\ldots,2,1,2,\ldots,d$} 
     \State Sample $f(\boldsymbol\Xi_{k})$, where $\boldsymbol\Xi_k$ is according to \eqref{eq:set-kron}. \Comment{$\mathcal{O}(n_{\xi}r^2)$ samples of $f$}
     \State Solve $F_{\neq k} f^{(k)} = f(\boldsymbol\Xi_{k})$ using the matrix \eqref{eq:frame-kron}.
     \State Compute $\mathcal{I}_{\ge k}$ from \emph{maxvol} on $(\overline{F}_{\ge k})^\top$ as defined in \eqref{eq:cross-iright}, and $\mathcal{I}_{\le k}$ from \emph{maxvol} on $\overline{F}_{\le k}$ as defined in \eqref{eq:cross-ileft}.
     \State Let $\Xi_{>k-1} = [\Xi^{(k)} \times \Xi_{>k}]|_{\mathcal{I}_{\ge k}}$, $\Xi_{<k+1} = [\Xi_{<k} \times \Xi^{(k)}]|_{\mathcal{I}_{\le k}}$, and $F_{>k-1} = \overline{F}_{\ge k}(:,\mathcal{I}_{\ge k})$, $F_{<k+1} = \overline{F}_{\le k}(\mathcal{I}_{\le k},:)$.
   \EndFor
 \EndWhile
\end{algorithmic}
\end{algorithm}

\section{Tensor Train algebra}\label{sec:ttalg}
Once a TT decomposition is constructed, simple operations can be performed directly with the TT cores.
Besides the fast integration \eqref{eq:int-rec},
additions of functions approximated by their TT formats,
pointwise products and actions of linear operators can be written in the TT format with explicitly computable TT cores, and hence a linear complexity in $d$ \cite{osel-tt-2011}.
For example, addition of functions $\tilde f(\xi)$ and $\tilde g(\xi)$ defined by their TT cores $\{F^{(k)}\}$ and $\{G^{(k)}\}$ and TT ranks $(r_0,\ldots,r_d)$ and $(p_0,\ldots,p_d)$, respectively, is realized by the TT decomposition
$$
\tilde f(\xi) + \tilde g(\xi) = \sum_{s_0,\ldots,s_d=1}^{r_0+p_0,\ldots,r_d+p_d} H^{(1)}_{s_0,s_1}(\xi^{(1)}) \cdots H^{(d)}_{s_{d-1},s_d}(\xi^{(d)}), \quad
H^{(k)}(\xi^{(k)}) = \begin{bmatrix}F^{(k)}(\xi^{(k)}) & 0 \\ 0 & G^{(k)}(\xi^{(k)})\end{bmatrix}.
$$
However, such explicit decompositions are likely to be redundant.
For example, we can immediately compress $H^{(1)}$ to
$[F^{(1)}(\xi^{(1)}) \quad G^{(1)}(\xi^{(1)})]$, and similarly the summation over $s_d$ collapses $H^{(d)}$.
In general, a quasi-optimal approximate \emph{re-compression} of a TT decomposition \cite{osel-tt-2011} can be performed in $\mathcal{O}(dn_{\xi} r^3)$ operations by using truncated Singular Value decompositions (SVD) in a recursive manner similar to \eqref{eq:int-rec}.

Linear operators acting on multivariate functions can also be decomposed (or approximated) in a similar TT format that enables a fast computation of their action.
An operator $A \in \mathcal{L}(\mathcal{F},\mathcal{F})$ with $\mathcal{F} = \mathcal{F}^{(1)} \times \cdots \times \mathcal{F}^{(d)}$
can be approximated by a TT operator $\tilde A$ of the form
\begin{equation}\label{eq:ttm}
 \tilde A = \sum_{t_0,\ldots,t_d=1}^{R_0,\ldots,R_d} A^{(1)}_{t_0,t_1} \otimes A^{(2)}_{t_1,t_2} \otimes \cdots \otimes A^{(d)}_{t_{d-1},t_d},
\end{equation}
where $A^{(k)}_{t_{k-1},t_k} \in \mathcal{L}(\mathcal{F}^{(k)},\mathcal{F}^{(k)})$ is an operator acting on $F^{(k)}(\xi^{(k)})$.
The image $\tilde A \tilde f$ can now be written as a TT decomposition
$$
(\tilde A \tilde f)(\xi) = \sum_{m_0,\ldots,m_d=1}^{R_0r_0,\ldots,R_dr_d} B^{(1)}_{m_0,m_1}(\xi^{(1)}) \cdots B^{(d)}_{m_{d-1},m_d}(\xi^{(d)}), \qquad B^{(k)}(\xi^{(k)}) = \left[(A^{(k)}_{t_{k-1},t_k} F^{(k)}_{s_{k-1},s_k})(\xi^{(k)})\right],
$$
where $t_k = \mathrm{mod}(m_k-1,R_k)+1$, and $s_k = \lfloor (m_k-1)/R_k \rfloor + 1$.
Linear equations can be solved by the Alternating Least Squares.
We iterate over $k=1,\ldots,d$, solving in each step $\tilde A \tilde f = \tilde b$ with respect to the TT core $\mathbf{F}^{(k)}$ in the representation of $\tilde f$.
Constructing a vector-function $f^{(\neq k)}: \mathbb{R}^d \rightarrow \mathbb{R}^{r_{k-1}n_{\xi}r_k}$ with elements
\begin{equation}\label{eq:als-frame}
 f^{(\neq k)}_{s_{k-1} i s_k}(\xi) = \sum_{s_0,\ldots,s_{k-2}} F^{(1)}_{s_0,s_1}(\xi^{(1)}) \cdots F^{(k-1)}_{s_{k-2},s_{k-1}}(\xi^{(k-1)}) \cdot \ell_i(\xi^{(k)}) \cdot \sum_{s_{k+1},\ldots,s_d} F^{(k+1)}_{s_{k},s_{k+1}}(\xi^{(k+1)}) \cdots F^{(d)}_{s_{d-1},s_d}(\xi^{(d)})
\end{equation}
similarly to \eqref{eq:cross-frame}, we can write $\tilde f(\xi) = f^{(\neq k)}(\xi) f^{(k)}$,
and solve a linear system
$$
A_k f^{(k)} = b_k, \quad \mbox{where} \quad A_k = \left[\langle f^{(\neq k)}_s, \tilde A f^{(\neq k)}_t \rangle\right], \quad b_k = \left[\langle f^{(\neq k)}_s, \tilde b \rangle\right]
$$
in each step.
If $\tilde A$ and $\tilde b$ are given by TT decompositions,
$A_k$ and $b_k$ can be computed recurrently in the course of iteration using matrix products similar to \eqref{eq:cross-ileft}, \eqref{eq:cross-iright}.
One full iteration over $k=1,\ldots,d$ needs $\mathcal{O}(dn_{\xi}^2 R^2 r^2 + d n_{\xi} R r^3)$ operations \cite{holtz-ALS-DMRG-2012}, where $R:=\max_k R_k$.

\end{document}